\documentclass[9pt,amstex]{article}
\usepackage{amssymb,amsmath} 
\usepackage{latexsym}
\usepackage{url}
\usepackage{hyperref}
\usepackage[usenames,dvipsnames]{color}

\newtheorem{dfn}{Definition}[section] 
\newtheorem{rmk}{Remark}[section]
\newtheorem{rmks}{Remarks}[section]
\newtheorem{thm}{Theorem}[section] 
\newtheorem{cor}{Corollary}[section]
\newtheorem{prop}{Proposition}[section] 
\newtheorem{lem}{Lemma}[section]
\newtheorem{exs}{Examples}[section]
\newtheorem{ex}{Example}[section]

\def\cyclic{\mathop{\kern0.9ex{{+}
\kern-2.2ex\raise-.28ex\hbox{\Large\hbox
{$\circlearrowright$}}}}\limits}

\def\buildrel#1_#2^#3{\mathrel{\mathop{\kern 0pt#1}\limits_{#2}^{#3}}}

\newcommand{\Pf}{{\em Proof}. }
\newcommand{\EPf}
{%
\mbox{}%
\nolinebreak%
\hfill%
\rule{2mm}{2mm}%
\medbreak%
\par%
}

\newcommand{\End}{\mbox{$\mathtt{End}$}}

\newcommand{\id}{\mbox{$\mathtt{id}$}}
\newcommand{\Aff}{\mbox{Aff}}

\newcommand{\Ad}{\mbox{$\mathtt{Ad}$}}
\newcommand{\pr}{\mbox{$\mathtt{pr}$}}
\newcommand{\ad}{\mbox{$\mathtt{ad}$}}
\newcommand{\rad}{\mbox{$\mathtt{rad}$}}

\newcommand{\C}{\mathbb C}

\renewcommand{\L}{\mathbb L}

\newcommand{\D}{\mathbb D}

\newcommand{\fM}{{\mathfrak M} {}}

\newcommand{\R}{\mathbb R}
\newcommand{\K}{\mathbb K}

\renewcommand{\sp}{{\mathfrak{sp}}} 
\newcommand{\g}{{\mathfrak{g}}{}}

\newcommand{\fZ}{{\mathfrak{Z}}{}}

\renewcommand{\k}{{\mathfrak{k}}{}}

\newcommand{\fL}{{\mathfrak{L}}{}}

\newcommand{\fG}{{\mathfrak{G}}{}}

\newcommand{\n}{{\mathfrak{n}}{}} 
 
\newcommand{\s}{{\mathfrak{s}}{}}

\newcommand{\h}{{\mathfrak{h}}{}} 
\newcommand{\fK}{{\mathfrak{K}}{}}

\newcommand{\CO}{{\cal O}{}}
\newcommand{\CN}{{\cal N}{}}
\newcommand{\CG}{\g{}}

\newcommand{\CM}{{\cal M}{}}

\newcommand{\CP}{\mathcal P}
\newcommand{\CS}{\mathcal S}

\newcommand{\CK}{\mathcal K}
\newcommand{\CL}{{\cal L}{}} 
\newcommand{\CR}{{\cal R}{}} 

\newcommand{\CJ}{{\mathcal J}}
\newcommand{\CI}{\mathcal I}

\newcommand{\CC}{\mathcal C}

\newcommand{\CZ}{\mathcal Z}

\newcommand{\ddto}{\left.\frac{{\rm d}}{{\rm d}t}\right|_0}

\def\cref#1{Corollary~\ref{#1}}

\title{Kinematical Lie algebras and symplectic symmetric spaces I\\ Lie algebraic aspects}

\author{
    Pierre Bieliavsky\thanks{Institut de Recherche en Math\'ematique et Physique, Universit\'e Catholique de Louvain, Chemin du Cyclotron, 2, 1348 
    Louvain-la-Neuve, Belgium. Email: Pierre.Bieliavsky@uclouvain.be} \and
    Nicolas Boulanger\thanks{Service de Physique de l'Univers, Champs et Gravitation, Université de Mons -- UMONS, 20 place du Parc, 7000 Mons, Belgium. Email: nicolas.boulanger@umons.ac.be}
}
\date{}
\begin{document}
\maketitle

\begin{abstract}

The aim of this note is to present a close relation between kinematical 
Lie algebras and symmetric spaces in a symplectic context: to every kinematical 
Lie algebra is canonically associated a symplectic symmetric space. 
For non-flat symmetric spaces, 
this correspondence is one-to-one onto a specific class of symplectic 
symmetric spaces whose structure we describe in details. 
In particular, the transvection Lie algebra of such a symmetric space is either 
three-graded or of the Poincar\'e type. The denomination ``Poincar\'e type'' 
refers to symplectic symmetric spaces characterized by a property that generalizes 
the fact that the classical Poincar\'e group $SO_o(1,D)\ltimes\R^{D+1}$ turns out 
to be the transvection group of an unexpected purely symplectic symmetric 
space structure on the cotangent bundle of the hyperbolic space $SO_o(1,D)/SO(D)$ 
(i.e., the mass-shell orbit).
The Lie triple system associated with every such symplectic symmetric space 
is of Jordan type (in the sense of W. Bertram), i.e., it is a homotope of the 
Lie triple system associated with a Hermitian symmetric space. 
However, the class of those Jordan-Lie triple systems associated with 
a classical kinematical Lie algebra is not stable under the natural operations 
of homotopies and dualities defined on Jordan-Lie triple systems. 
In order to restore stability, we need to introduce a natural generalization 
of the notion of kinematical Lie algebras, which is the framework where the present 
work is formulated. The last section of this work presents some remarks on the 
coadjoint orbits that are naturally associated with symplectic symmetric spaces.

\end{abstract}

\section{Introduction}

\subsection{Physical aspects}

Of central importance in physics is the notion of (equivalence class of) 
inertial reference systems, or inertial frames. To some extent, the 
set of reference frames together with the group of transformations 
relating them determine the global geometry of spacetime. The laws of physics 
are invariant under the transformations relating the inertial frames of 
a given equivalence class of such frames. 
In the Newtonian formulation of mechanics, two inertial frames 
are related by a transformation of the Galilean group. 
With the works of Lorentz, Poincar\'e and Einstein, it was 
understood that two inertial frames can more accurately be related by 
what is now called a Lorentz, or Poincar\'e, transformation, 
a transformation which leaves the speed of light invariant.
The Lie algebras of these transformation groups are maximal, 
in the sense that, in a space-time with three space dimensions and one time 
direction, they contain three generators of ``translation", 
three generators of ``boost", three rotation generators that form 
a basis of $\mathfrak{s}\cong \mathfrak{so(3)}$,
and one time-translation generator. Altogether, this gives ten generators,
ten being the number of linearly independent Killing vectors of the maximally symmetric 
(pseudo)Riemannian spaces of dimension four \cite{eisenhart1966riemannian}.

In two seminal papers, Bacry, L\'evy-Leblond and Nuyts 
\cite{Bacry:1968zf,Bacry:1986pm}
determined the various possible relativity principles that 
a four-dimensional, isotropic spacetime could accommodate, thereby 
defining the various possible classes of inertial frames. 
In other words, they asked whether there could exist theoretically possible 
relativity principles that differ from the Galilean or the Lorentzian ones 
discussed above. The various possible transformation groups  
expressing the equivalence relations among inertial frames within each 
possible class were qualified as \emph{kinematical} in 
\cite{Bacry:1968zf,Bacry:1986pm}, and the main 
result of these two papers was the classification of the 
kinematical algebras pertaining to four-dimensional, isotropic 
spacetimes.
It was also observed in \cite{Bacry:1968zf,Bacry:1986pm} that, 
to every kinematical Lie algebra $\mathfrak{g}$ they classified, 
there always exists a six-dimensional subalgebra 
$\tilde{\mathfrak{h}}\subset \mathfrak{g}\,$, suggesting the 
existence of a four-dimensional homogeneous space 
$G/\tilde{H}$ describing a possible isotropic spacetime, 
the rotation algebra $\mathfrak{so(3)}$ 
always being contained in $\tilde{\mathfrak{h}}$.
As a vector space, any kinematical Lie algebra $\mathfrak{g}$ 
appearing in the classification of \cite{Bacry:1968zf,Bacry:1986pm}
can be decomposed as 
$\mathfrak{g}=\mathfrak{s}\oplus \mathcal{P}\oplus \mathcal{Z}$,
where $\mathfrak{s}\cong \mathfrak{so}(3)$, and  
$\mathcal{P}=
\mathcal{P}_0\oplus \mathcal{P}_1$. 
The subspace $\mathcal{Z}$ is a line, while both $\mathcal{P}_0$ 
and $\mathcal{P}_1$ are $\mathfrak{so}(3)$ modules in 
the vector representation.

The classification of Bacry, Lévy-Leblond, and Nuyts was reconsidered 
later \cite{ofarrill2017defo,ofarrill2017higher,andrzejewski2018} 
in higher and lower dimensional space-times, replacing the isotropy 
algebra $\mathfrak{so}(3)$ by $\mathfrak{so}(D)$, with 
the number of space dimensions $D\geq 1$ and the dimension of 
$\CP_0$ and $\CP_1$ being equal to $D$, for the vector representation 
of $\mathfrak{so}(D)$. \textcolor{blue}{The detailed relations 
between the various Klein pairs $( \mathfrak{g}, \tilde{\mathfrak{h}})$ 
and the various possible isotropic spacetimes of the form $G/\tilde{H}$ 
were given in \cite{Figueroa-OFarrill:2018ilb}. In the latter reference
it was stressed  that not all Klein pairs effectively give rise to homogeneous 
spaces, neither is it necessary that they should be geometrically realisable. 
See also \cite{Figueroa-OFarrill:2017sfs,Figueroa-OFarrill:2022nui,Morand:2023emw} 
for a summary and related results.}

\textcolor{blue}{In the present paper, we show that, already for the cases 
$\mathfrak{so}(D)$ with $D>3\,$, there is a symplectic symmetric 
space \cite{BiThese,BCG2} that governs the structure of any of these classical 
kinematical Lie algebras. Besides, it appears
necessary to further extend the notion of kinematical Lie algebra to what we call
\emph{generic} kinematical Lie algebra that encompass the cases where 
$\mathfrak{s}$ can be any Lie algebra, not necessarily $\mathfrak{so}(D)$.
The reason why such an extension should necessarily appear is spelled out in 
next section \ref{subsec:Math}. The notion of generic kinematical Lie algebra 
defines remarkable families of Lie algebras and symmetric spaces 
that persist in dimensions $D\leqslant 3$, and that were specifically 
considered by Figueroa-O'Farrill and collaborators  
in the classical case $\mathfrak{s}\cong \mathfrak{so}(D)$.}

The classification of the kinematical algebras is not a purely academic question. 
Although Minkowski's spacetime is a well adapted arena where  to describe 
physical processes where particles can travel at velocities approaching 
the speed of light, the Newtonian space-time is perfectly adapted 
(and used every day) to host physical phenomena where the velocities 
of particles are negligible compared to the speed of light.
On the other hand, modern cosmology (see e.g. \cite{weinberg2008cosmology} 
for a textbook) indicates that 
the geometry of our early, inflationary universe can be approximated 
by a de-sitterian spacetime, with a positive, constant curvature. 
Finally, in the context of string theory, J. Maldacena \cite{maldacena1998large}
formulated a conjecture that gives a prime importance 
to anti-de sitterian (AdS) spacetimes of various dimensions.
From a practical point of view, depending on the 
nature of the physical system under study, one given structure of spacetime 
will be preferred to another.

\subsection{Mathematical aspects}
\label{subsec:Math}

\subsubsection{Symplectic symmetric spaces}
\noindent Symplectic symmetric spaces were introduced by one of us et al. in the mid 
nineties \cite{BiThese,BCG2}. The idea is
to define a class of symplectic manifolds admitting a \emph{preferred} symplectic 
linear connection. There are many ways to do this, e.g. 
requiring a compatibility with a (pseudo) Riemannian metric, considering critical 
solutions of a variational problem on symplectic connections, etc.
(see e.g. \cite{BCGS} and references there           in).
Imposing the symplectic manifold to admit a large \cite{BCG2}
set of symmetries is one of them. 
This choice was originally motivated, on the first hand,  by the positive solution 
due to Sekigawa and Van Hecke \cite{SVH} of a conjecture of S. Kobayashi asserting 
that a compact K\"ahler manifold whose (local) geodesic symmetries 
(w.r.t. the K\"ahler metric) are symplectic (w.r.t. the K\"ahler 2-form) must 
be a (local) Hermitian symmetric space. And on the second hand by a conjecture 
due to A. Weinstein in the context of invariant star-product on Hermitian symmetric 
spaces \cite{W}.

\noindent A symplectic symmetric space is a Fedosov manifold, i.e., a symplectic 
manifold $(M,\omega)$ equipped with a torsionfree symplectic linear connection 
$\nabla$ (i.e. $\nabla\omega=0$) enjoying the property that every $\nabla$-geodesic 
symmetry centred at any point of $M$ is well-defined as a global affine transformation 
of $(M,\nabla)$. Locally, this last property amounts to saying that the curvature 
$(3,1)$-tensor $R^{\nabla}$ of the torsionfree connection $\nabla$ is parallel: 
$\nabla R^{\nabla}=0\,$; beware the fact that this last equation is not equivalent 
nor implied by a Bianchi identity. Affine symmetric spaces can be described and 
studied within several (equivalent) geometrical contexts depending on the user's 
needs and preferences. There are three main approaches. 
The first one is ``differential geometrical'', i.e., based on the datum of the 
connection. The second one is Lie theoretical, i.e., where the symmetric manifold 
is considered as a homogeneous space $G/H$, and the third one, 
adopted in the present article, consists in O. Loos' approach to symmetric 
spaces within the Jordan algebraic context \cite{Lo}.
In this last framework, a symmetric space appears as a manifold with 
multiplication, generalising the notion of Lie group; see Remark \ref{RMKLG} below.

\subsubsection{Semisimple symplectic symmetric spaces, their homotopes and dualities}
\noindent In contrast with Riemannian symmetric spaces (i.e., Riemannian 
manifolds whose centered geodesic symmetries are globally defined and all isometrical),
symplectic symmetric spaces are generally not homogeneous spaces of semisimple Lie groups. 
And conversely, a simple Lie group does generally not cover any symplectic symmetric 
space; for the structure and classification of semisimple symplectic symmetric spaces, 
see \cite{Bi98bis}. 
However, those that are semisimple constitute an important class of symplectic 
symmetric spaces. A first reason being that the class of semisimple symplectic 
symmetric spaces is acted on by several dualities. 
The compact-noncompact duality for Hermitian symmetric spaces is one of them 
(see e.g. \cite{Helgason1978}). 
Another one, called Hermitian-Cayley type, is the 
duality between Hermitian symmetric spaces and the 
causal symmetric spaces of Cayley type \cite{faraut1994analysis,kaneyuki1985}. 
To some extent, one can consider the class of hyper-K\"ahler symmetric spaces 
as the fixed points of those dualities.
A second reason of importance of the semisimple class is that it naturally 
lies in the framework of Jordan algebra theory, which allows for considering 
\emph{homotopies} between spaces as special kinds of deformations within 
the specific class of the \emph{Jordan-Lie triple systems}. 
Roughly speaking, given a simple symplectic symmetric space, 
the set of its homotopes naturally consists in an algebraic variety \cite{Bertram,BeBi1}, 
called in Bertram's terminology the \emph{structure variety}. 
The interior points of such a structure variety are simple symmetric spaces, 
while the ones on the boundary not anymore (see \cite{BeBi1, BeBi2} for structure 
and full classification). A continuous path in the structure variety 
from an interior point to a boundary point can be thought of as some 
kind of generalized Inonu-Wigner contraction. 

\noindent Nevertheless, the above mentioned dualities do not always correspond to 
geometric constructions. An example of a duality which can geometrically be 
realised is the one of compact-noncompact Hermitian symmetric spaces. 
For those, the duality is based on holomorphic holonomy-equivariant embeddings 
generalizing the $U(1)$-equivariant holomorphic embedding of the hyperbolic 
plane into the Riemann sphere.
The Hermitian-Cayley type duality is not as clearly associated to 
such a geometrical construction. All this brings us to kinematical Lie algebras.

\subsubsection{Kinematical Lie algebras and the Hermitian-Cayley type duality}
\noindent As a byproduct of the main result of the present work, a generic 
kinematical Lie algebra corresponds to a homotope of a Hermitian symmetric space. 
By \emph{generic} we mean the following.  In space-time dimension greater than
or equal to five, the fact that the isotypical component of the natural 
representation of the rotation Lie algebra in its anti-symmetric square is 
empty separates the isomorphism classes of kinematical Lie algebras into 
distinguished families that will described below. 
By \emph{generic}, we mean a kinematical Lie algebra 
belonging to one of those families. 
In smaller space-time dimensions, the aforementioned generic families subsist, 
but there are also other types \cite{Figueroa-OFarrill:2018ilb}, 
due to the fact that the above property of the isotypical component does 
not hold anymore.
These extra, lower dimensional classes, are not considered in the present work.
The four-dimensional case will be reconsidered at the light of the present 
context in a forthcoming work. The complete classification of kinematical 
Lie algebras with isotropy algebra $\mathfrak{so}(D)$ with generic space dimension
$D$ was given in \cite{ofarrill2017defo,ofarrill2017higher,andrzejewski2018}.

\noindent This being said, the other elements and their natural Jordan-duals 
(in the sense of the previous subsection)
of the structure variety containing those Hermitian symmetric spaces 
associated with the usual kinematical Lie algebras 
are lost when restricting the Levi factor of the holonomy 
Lie algebra to be isomorphic to the rotation Lie algebra $\mathfrak{so}(D)$. 

\noindent For this reason: in order to define a class of symmetric spaces 
containing the above-mentioned kinematical Hermitian symmetric spaces
and which is stable under Jordan-Lie homotopy and symmetric space dualities, 
we need to relax the condition on the Levi factor of the holonomy to be 
isomorphic to the rotations by only requiring the emptiness condition spelled out 
above on the isotypical component of the natural representation of the 
Lie algebra that replaces the rotation algebra of the usual case.
Precisely, we formulate our

\begin{dfn}\label{KLA} A {\bf generic\footnote{the adjective ``generic'' refers to two things. First, to the fact that the dimension $D$ of $V$ is generic e.g. greater or equal to four in the classical case where $\s$ is isomorphic to the rotation Lie algebra $\mathfrak{so}(D)$.
And second to the fact that the Lie algebra $\s$ is not defined as being isomorphic to the rotation Lie algebra.}  kinematical Lie algebra} is a triple 
$(\g,\s,V)$ where 
\begin{enumerate}
\item $\g$ is a finite-dimensional Lie algebra over $\K=\R$ or $\C$,
\item $\s$ is a  Lie subalgebra of $\g$, and 
\item $V$ is a faithful absolutely simple\footnote{Recall that a simple module 
$M$ over an 
algebra $A$ on a ground field $\K$ is called absolutely simple if 
the only $A$-commuting $\K$-endomorphisms of $M$ are in the ground field $\K$: 
$\End_\K(M)/\End_A(M)\simeq \K$.}  
finite-dimensional $\s$-module over $\K$ such that 
\begin{enumerate}
\item the \emph{weak isotypical component} (see below) of $V$ in 
$\Lambda^{2}(V)$ is empty.
\item The action of $\s$ on $V$ preserves a non-degenerate scalar product (not necessarily positive definite) on $V$.
\end{enumerate}
\end{enumerate}
The triple is moreover conditioned by the requirement that 
$\g$ admits a vector space decomposition
\begin{equation}\label{KIN}
\g\;=\;\CZ\,\oplus\,\s\,\oplus\,\CP
\end{equation}
where 
\begin{enumerate}
\item[4.]
$\CP$ is stable under $\s$ and isomorphic to  the reducible $\s$-module 
$V\oplus V$, and 
\item[5.] $\CZ$ is a line which centralizes $\s$ in $\g$.
\end{enumerate}
\end{dfn}

\noindent By \emph{weak isotypical component} we mean
\begin{dfn}\label{weak}
\noindent Let $W$ be  a (not necessarily semisimple) $\s$-module and $L$ an 
irreducible $\s$-module. The {\bf weak isotypical component} $W_{(L)}$ of $L$ 
in $W$ is defined as  the vector sum in $W$ of all sub-modules of $W$ that are 
isomorphic to $L$. 
\end{dfn}

\noindent In the sequel, we will abusively use the adjective ``isotypical'' 
for ``weakly isotypical''.
\begin{rmks}
\noindent 1. In the context of Definition \ref{weak},  
the $\s$-module $W_{(L)}$ does not necessarily admit a supplementary 
$\s$-submodule in $W$.

\noindent 2. The requirement on $V$ to be absolutely simple 
(as opposed to only simple) is not essential but it simplifies the discussion. 
\end{rmks}

\noindent Although the deformational aspects of the present study will 
be deferred to a further article, we anticipate that the main result
of the present work enables to tie the notion of (necessarily generalized) 
kinematical Lie algebra to a specific class of symplectic symmetric spaces 
that is stable under homotopy and duality in the context of Jordan triple systems. 
We will precise this in a forthcoming work.

\subsection{Structure of the present work}

\noindent In its essence, this article is interdisciplinary as it realizes 
a bridge between a theoretical physics notion and a differential geometric one.  
In order to reach a wider audience, in particular our colleagues in the 
theoretical physics community, we decide to write, in a first part of this article, 
a concise and essentially self-contained introduction to symmetric spaces and their 
symplectic version.
In particular, up to basics in differential geometry, we made an effort to present 
proofs of the results that are relevant to this work.

\subsubsection{Basics on symmetric spaces and Symplectic symmetric spaces}

Section \ref{SS} consists in an introduction to symmetric 
spaces and their symplectic analogues.
All the material in this introduction has been well established for several decades, 
the main references being \cite{KN1,KN2}, \cite{Lo}, \cite{Helgason1978}, 
\cite{BiThese}, \cite{Bi2} and \cite{BCG2}.

\noindent We start by the notion of symmetric space within the  
approach of O. Loos. We then explain the correspondence between symmetric spaces
and their tangent analogues: involutive Lie algebras (iLa's). 
This correspondence generalizes to symmetric space Lie's third theorem.

\noindent We then proceed by introducing the main actor of this work, the notion of 
symplectic symmetric space, and describe, in the same lines of ideas as in 
the first step, their tangent analogues: symplectic iLa's (siLa's). 
We detail the correspondence (analogous to Lie's third theorem) that integrates 
siLa's to symplectic symmetric spaces. 
We also recall how a symplectic symmetric space can, up to automorphism, 
uniquely be decomposed into elementary pieces. This result is analogue 
to the de Rham decomposition theorem in Riemannian geometry\cite{de1952reductibilite}.

\subsubsection{Generic kinematical Lie algebras correspond 
to symplectic symmetric spaces}

\noindent Section \ref{HEART} constitutes the heart of the present work: 
we first show that, to every generic kinematical 
Lie algebra $(\g,\s,V)$ as defined in Definition \ref{KLA} is canonically 
associated a simply connected  symplectic symmetric space $(M,\omega,\nabla)$ of 
dimension $\dim\CP$ on which the simply connected Lie group $G$ admitting $\g$ 
as Lie algebra acts by automorphisms, i.e., by symplectic diffeomorphisms of 
$(M,\omega)$ that preserve the connection $\nabla$.
As a by-product of the construction we prove that every generic kinematical 
Lie algebra $(\g,\s,V)$ satisfies $[\CZ,\CP]\subset\CP$, i.e., $\CZ$ acts on $\CP$.

\noindent We then concentrate on describing the fine structure of the 
symplectic symmetric spaces associated with our generic kinematical Lie algebras. 
By what was presented in Section \ref{SS}, it is sufficient to consider 
the generic kinematical Lie algebras for which the associated simply connected 
symplectic symmetric space is indecomposable and non-flat (i.e., $R^\nabla\neq0$), 
which we assume from now on.

\noindent For those, it turns out that the action of $\CZ$ on $\CP$ is 
either (complex) semisimple or nilpotent (no mixing).
We entirely describe the structure in the case the action of $\CZ$ on $\CP$ is  
(complex) semisimple.
For instance, those generic kinematical Lie algebras that are such that $\s$ is 
semisimple are \emph{exactly} the simple transvection Lie algebras of simple symplectic 
symmetric spaces. These spaces are classified in \cite{Bi98bis} and there structure is 
described: they are either Hermitian symmetric spaces, para-Hermitian causal symmetric 
spaces of Cayley type, or hyper-K\"ahler symmetric spaces.

\noindent We describe the spaces for which the action of $\CZ$ on $\CP$ 
is nilpotent under the condition that $\s$ is contained in a Levi factor of 
$\g$; this ``Levi'' condition for instance holds when $\s$ is compact semisimple. 
A simply connected symplectic symmetric space associated with such a generalized 
kinematical Lie algebra is called of the \emph{Poincar\'e type}. 
We prove that every such symplectic symmetric space is equivariantly symplectomorphic 
to the cotangent bundle of a simple symmetric space. 
This result classifies those spaces. In particular, when $\s$ is compact: 
the Poincar\'e types are exactly the cotangent bundles of Cartan's Riemannian 
symmetric spaces. Combining these results with the $\CZ$-semisimple action yields 
a complete classification in all the cases where $\s$ is semisimple.

\subsubsection{Hamiltonian aspects}

\noindent Section \ref{HAM} is two-fold: it first provides a concise presentation 
of standard basic facts on homogeneous symplectic spaces, and, secondly, 
it presents some features of the dynamical aspects of our present work. 
The original results of this section consist in two facts. First, the action of 
the group (necessarily Lie) generated by the geodesic symmetries on the 
(non-flat, indecomposable) symplectic symmetric space associated with a generalized 
kinematical Lie algebra is strongly Hamiltonian in the sense that it defines an 
equivariant moment map.
And second, the Lie algebra of this Lie group also underlies a structure of 
generic kinematical Lie algebra.

\subsubsection{Conclusions and perspectives}

\noindent In the last Section \ref{CONCL}, 
we summarise what has been done in the present article and suggest some of its 
developments.

\vspace{3mm}

\noindent {\bf Acknowledgement.} We thank the referee for having 
pointed out several inconsistencies in the first version of the manuscript 
as well as for having suggested substantial improvements of the presentation. 

\section{Basics on symmetric spaces and symplectic symmetric spaces}\label{SS}
\subsection{Symmetric spaces}
\noindent We first recall some basic facts about (affine) symmetric spaces. The 
references for the results appearing in this section are 
\cite{Lo}, \cite{KN1}, \cite{BiThese} and \cite{Bi2}. 
Roughly speaking, a (locally) symmetric space is an affine manifold 
whose geodesic symmetries are affine transformations. More precisely:
\begin{dfn}\label{SSDEF}\cite{Lo} A {\bf symmetric space} 
is a pair $(M,s)$, where $M$ is a smooth connected 
manifold, and where $s : M \times M \to M$ is a smooth map such 
that 
  
\begin{enumerate} 
\item[(i)] for all $x$ in $M$, the partial map $s_x : M \to M : y \mapsto 
s_x (y) := s(x,y)$ is an involutive diffeomorphism of $M$ 
called the {\bf symmetry} at $x$. 
\item[(ii)] For all $x$ in $M$, $x$ is an isolated fixed point of $s_x$. 
\item[(iii)] For all $x$ and $y$ in $M$, one has $s_xs_ys_x=s_{s_x (y)}$ ({\bf 
Fundamental Jordan identitity}). 
\end{enumerate}
\end{dfn}

\begin{dfn} Two symmetric spaces  $(M,s)$ 
and $(M',s')$ are {\bf isomorphic} if there exists a 
diffeomorphism $\varphi: M \rightarrow M'$ such that 
$\varphi s_x=s'_{\varphi (x)} \varphi$. Such a $\varphi $ is called an {\bf
isomorphism} of $(M,s)$ onto $(M',s')$. When 
$(M,s) = (M',s')$, one talks about {\bf automorphisms}. 
The group of all automorphisms of the  symmetric space 
$(M,s)$ is denoted by $Aut(M,s)$.
\end{dfn}

\begin{ex}\label{RMKLG}
An important class of symmetric spaces is the one of Lie groups. Given a Lie 
group $\fG$ with group law
$(x,y)\mapsto xy$, the map
$$
s:\fG\times\fG\to\fG:(x,y)\mapsto xy^{-1}x
$$
defines a structure of symmetric space on the group manifold $\fG$.
\end{ex}

\begin{prop} On a symmetric space $(M,s)$, there exists 
one and only one affine connection $\nabla$ which is invariant under the 
symmetries. 

\noindent It is explicitly given by the following formula :
\begin{equation}\label{CONNLOOS}
\left(\nabla _X Y\right)_{x}\;=\;\frac{1}{2}
\left[X\,,\,Y+s_{x_{\star}}Y\right]_{x}
\end{equation}
at every point $x$ of $M$ and for all  tangent smooth vector fields $X$ and $Y$ on $M$.

\vspace{2mm}

\noindent The linear connection $\nabla$ enjoys the properties of being 
torsionfree and such that its curvature 
$(3,1)$-tensor $R^{\nabla}$ is parallel.  

\noindent Moreover, the automorphism group of the symmetric space $(M,s)$ 
coincides with the group of the affine transformations of the connection $\nabla$: 
$$
\mbox{Aut}(M,s)\;=\;\mbox{Aff}(M,\nabla)\;.
$$
In particular, the automorphism group $\mbox{Aut}(M,s)$ is a Lie group of 
transformations of $M$.
\end{prop}

\noindent Using basic facts in differential geometry (see e.g. \cite{KN1}), 
the proof of the above proposition is straightforward: it suffices to check that formula 
(\ref{CONNLOOS}) indeed defines a covariant derivative 
in the tangent bundle of $T(M)$, which is a routine computation.
One  readily verifies the announced properties using the explicit formula 
(see \cite{BiThese} and \cite{bertelson2011affine}).

\begin{rmk}
\noindent The existence and uniqueness of a symmetry-invariant affine connection 
on every symmetric space was proved by Loos in \cite{Lo} volume I. For this reason, 
this canonical connection on such a symmetric space will be hereafter called the 
{\bf Loos connection}. Note, however, that the explicit formula for it was first 
given in \cite{bertelson2011affine}.
\end{rmk}

\begin{ex}
\noindent In the case of the flat symmetric space $\R^n$ equipped with its euclidean 
affine symmetries, the affine group is $\Aff(\R^n)\;=\;\mbox{GL}_n(\R)\ltimes\R^n$. One 
observe that in this case, the subgroup $\R^n$ also transitively acts on the space. And 
moreover, it is characteristic of the flat connection in the sense that there is only one 
connection on $\R^n$ that is invariant under the translation group $\R^n$. 
\end{ex}
\noindent This observation is actually general. First observe
\begin{prop}\label{AUTO}
The automorphism group $\mbox{Aut}(M,s)$ transitively acts on $M$.
\end{prop}
\Pf Consider any two points $x$ and $y$ in $M$. By connectedness, there exists a 
continuous path $\gamma:[0,1]\to M$ joining them:
$\gamma(0)=x$ and $\gamma(1)=y$. For every point $z$ of the curve $C:=\gamma([0,1])$, 
choose a $\nabla$-normal neighbourhood $U_{z}$ and 
extract (by compactness of $C$) a finite open cover 
$\{U_{i}:=U_{\gamma(t_{i})}\}_{t_{1}=0\,<t_{2}<...<t_{N}=1}$ of $C$ 
in $\cup_{z\in C}U_{z}$. In every intersection $U_{i}\cap U_{i+1}$ choose a 
point $p_{i}$ and consider geodesic arcs $\eta_{\gamma(t_{i})}^{p_{i}}$ joining 
$\gamma(t_{i})$ to $p_{i}$ and $\eta_{p_{i}}^{\gamma(t_{i+1})}$ joining $p_{i}$ 
to $\gamma(t_{i+1})$. The union of all these arcs then defines a broken geodesic 
line joining $x$ to $y$. Denote by $m^{+}_{i}$ the mid-point of the arc 
$\eta_{\gamma(t_{i})}^{p_{i}}$ and by $m^{-}_{i}$ the mid-point of the arc 
$\eta_{p_{i}}^{\gamma(t_{i+1})}$. Then, the affine transformation 
$s_{m^{-}_{N}}\circ s_{m^{+}_{N}}...s_{m^{-}_{1}}\circ s_{m^{+}_{1}}$ sends 
$x$ onto $y$.
\EPf

\noindent In a second step, we introduce the analogue of the translation 
group in the flat case:

\begin{dfn}
The {\bf transvection group} (or {\bf displacement group}) $G(M,s)$ of a 
symmetric space $(M,s)$  is defined as the subgroup of the automorphism group 
$Aut(M,s)$ generated by 
$\{ s_x \circ s_y \, ; \, x,y \in M \}$. 
\end{dfn}

\noindent As a byproduct of Proposition \ref{AUTO} and its proof, one then gets 
(see  \cite{Lo} vol. I pp. 88 Thm 2.8.):
\begin{thm}\label{TG}
\noindent The transvection group $G(M,s)$ of a symmetric space is a Lie group of 
transformations of $M$. It is characterized by the property to be the smallest 
subgroup of the automorphism group that acts transitively over $M$
and that is stable under the conjugation by a (and therefore all) symmetry $s_{o}$ 
($o\in M$) in $\mbox{Aut}(M,s)$: 
$$
\tilde{\sigma}:\mbox{Aut}(M,s)\to\mbox{Aut}(M,s):g\mapsto s_{o}\circ g\circ s_{o}\;.
$$
\end{thm}
\begin{rmks}
\begin{enumerate}
\item[(i)] The notion of \emph{transvection} generalizes the one of translation: 
in euclidean space, the composition of two centred symmetries is a translation.
\item[(ii)] Observe that for every transvection $g\in G(M,s)$ and every point $x\in M$, 
one has $gs_{x}g^{-1}=s_{g(x)}$. Indeed, at the level of generators, the fundamental 
Jordan identity (Definition \ref{SSDEF} item (iii)) implies 
$s_{y}s_{z}s_{x}s_{z}s_{y}=s_{y}s_{s_{z}x}s_{y}=s_{s_{y}s_{z}x}$.
\end{enumerate}
\end{rmks}

\begin{rmk}\label{REMPARAL} Using formula (\ref{CONNLOOS}) it is easy to check that a 
tensor field on a symmetric space which is invariant by the symmetries is necessarily 
parallel w.r.t. the Loos connection. 
\end{rmk}

\noindent As a conclusion of the present paragraph, we observe the canonical 
correspondence $(M,s)\mapsto G(M,s)$ between symmetric spaces and 
\emph{involutive Lie groups}, i.e., pairs $(G,\tilde{\sigma})$ where $G$ is a 
Lie group and where $\tilde{\sigma}$ is an involutive automorphism of $G$.

\noindent In the next section, we will see how from the above discussion emerges a 
generalization (c.f. Example \ref{RMKLG}) of Lie's third theorem for symmetric spaces.

\subsection{Involutive Lie algebras and de Rham's decompositions}

\begin{dfn}
An {\bf involutive Lie algebra} (abreviated ``iLa'') is a pair $(\g,\sigma)$ 
where $\g$ is a finite-dimensional real Lie algebra and where 
$$
\sigma:\g\to\g
$$
is an involutive ($\sigma^2=\id_\g$) automorphism of the Lie algebra $\g$.
\end{dfn}
\noindent Associated with $\sigma$, one has a vector decomposition of $\g$ into its 
$\pm1$-eigenspaces:
\begin{equation}\label{ILADECOMP}
\g\;=:\;\h\;\oplus\;\CP\quad\mbox{with}\quad\sigma\;=\;\id_\h\;\oplus(-\id_\CP)\;.
\end{equation}
One then observes
\begin{lem}
\noindent Given an iLa with associated decomposition (\ref{ILADECOMP}), 
one has the inclusions
$$
\begin{array}{c}
\left[\h\,,\,\h\right]\;\subset\;\h\;,\\
\left[\h\,,\,\CP\right]\;\subset\;\CP\;,\\
\left[\CP\,,\,\CP\right]\;\subset\;\h\;.
\end{array}
$$
Reciprocally, if a Lie algebra $\g$ admits a vector decomposition 
$\g\;=\;\h\;\oplus\;\CP$ satisfying the above inclusions, 
then it is underlain by an iLa structure $\sigma$ on $\g\,$.
\end{lem}
\begin{dfn}
Given two iLa's $(\g_i\,,\,\sigma_i)$ ($i=1,2$), a {\bf morphism} between them is a Lie 
algebra homomorphism
$$
\varphi:\g_1\to\g_2
$$
which intertwines the involutions:
$$
\varphi\;\circ\;\sigma_1\;=\;\sigma_2\;\circ\;\varphi\;.
$$
When $\varphi$ is bijective one refers to it as an {\bf isomorphism}.
\end{dfn}

\begin{dfn}\label{TRANSVECTIONLA}
\noindent An iLa $(\g,\sigma)$ is called a {\bf transvection iLa} if the additional 
conditions hold
\begin{enumerate}
\item $[\CP,\CP]\;=\;\h$ , and
\item the action of $\h$ on $\CP$ is faithfull.
\end{enumerate}
\end{dfn}
\noindent Every iLa is associated with a ``canonical'' transvection iLa: 
we have the following lemma whose proof is obvious.
\begin{lem}
Let $(\g,\sigma)$ be an iLa with decomposition $\g\;=\;\h\oplus\CP$. Then 
\begin{enumerate}
\item[(i)]
$$
\hat{\g}\;:=\;
[\CP,\CP]\;\oplus\,\CP\;,\quad\hat{\sigma}\;:=\;\id_{[\CP,\CP]}\;
\oplus\;\left(-\id_\CP\right)
$$
is a sub-iLa of $(\g,\sigma)$.
\item[(ii)] The centralizer $\n$ of $\CP$ in $[\CP,\CP]$ is an ideal of $\hat{\g}$ and 
the associated  exact sequence
$$
\n\;\longrightarrow\;\hat{\g}\;\stackrel{\underline{\pi}}
{\longrightarrow}\;\underline{\g}\;:=\;\hat{\g}/\n
$$
naturally induces on the quotient $\underline{\g}\;:=\;\hat{\g}/\n$ a structure of iLa.
\item[(iii)] The iLa defined in item (ii) is a transvection iLa.
\end{enumerate}
\end{lem}
\begin{dfn}
\noindent The holonomy of the iLa $(\g,\sigma)$ is defined as the Lie algebra 
$[\CP,\CP]/\n$. The iLa $(\g,\sigma)$ is called {\bf flat} when 
$[\CP\,,\,\CP]\subset\n$. Equivalently, it can be shown that the Loos 
connection of the corresponding symmetric space is flat.
\end{dfn}
\begin{dfn}
Given two iLa's $(\g_i\,\sigma_i)$ ($i=1,2$), their {\bf direct sum} is defined 
as the natural iLa structure on the direct sum of the Lie algebras $\g_1\oplus\g_2$, 
i.e., $(\g_1\oplus\g_2\,,\,\sigma_1\oplus\sigma_2)$. 

\noindent An iLa is called {\bf decomposable} if it is isomorphic to a direct 
sum of two non-trivial iLa's. It is called {\bf indecomposable} otherwise. 
\end{dfn}
It is not at all clear that two decompositions into indecomposable pieces of a 
given iLa are \emph{isomorphic}, i.e., that there exists an automorphism of the iLa 
swapping them. In the Riemannian case, it is a consequence of the de Rham decomposition 
theorem \cite{de1952reductibilite}. Happily, it holds in the pure affine case as well 
\cite{BiThese,BCG2}:
\begin{thm}\label{DERHAM}
Let $(\g,\sigma)$ be a transvection iLa. And let 
$$
(\g\,,\,\sigma)\;=\;\oplus_{i=1}^r(\g_i\,,\,\sigma_i)\;
=\;\oplus_{j=1}^{\overline{r}}(\overline{\g}_j\,,\,\overline{\sigma}_j)
$$
be two decompositions into indecomposable iLa's.
Then, $r\;=\;\overline{r}$ and there exist a permutation $\tau\in\mbox{Sym}(r)$ and an automorphism $\varphi$ of $(\g\,,\,\sigma)$ such that for all $i\in\{1,...,r\}$:
$$
\varphi(\g_i)\;=\;\overline{\g}_{\tau(i)}\;.
$$
\end{thm}
Now, as announced earlier, we have the
\begin{thm}\label{BIJSS}
The correspondence $(M,s)\mapsto G(M,s)$ induces a bijection between 
the isomorphism classes of simply connected symmetric spaces and the 
isomorphism classes of transvection involutive Lie algebras.
\end{thm}
\Pf
\noindent Associating a simply connected symplectic symmetric space 
to an iLa $(\g,\sigma)$, although elementary, fully uses standard 
techniques in differential geometry:
since the bijection (Lie's third theorem) that associates a connected simply 
connected Lie group $G$ to a Lie algebra $\g$ is functorial, to the pair 
$(\g,\sigma)$ is associated a pair $(G,\tilde{\sigma})$ with 
$\sigma:=\tilde{\sigma}_{\star e}$. The subgroup $G^{\tilde{\sigma}}$ of group
elements that are fixed by $\tilde{\sigma}$ is closed. Indeed it is the pre-image 
of the identity $e$ of $G$ under the continuous map 
$G\to G:g\mapsto g^{-1}\tilde{\sigma}(g)$. 
The subgroup $G^{\tilde{\sigma}}$ is therefore an embedded Lie subgroup of $G$ 
and so is its connected component $H:=G^{\tilde{\sigma}}_{0}$ containing the 
identity $e$. The set $G/H$ of left lateral classes (the ``coset'' space) of $H$ 
in $G$ therefore carries a unique structure of smooth manifold, 
homogeneous under the obvious action of $G$ on the classes:  
$g_{0}(gH):=g_{0}gH\,$. The first terms of the long exact sequence in homotopy 
$...\pi_{1}(G)\to\pi_{1}(G/H)\stackrel{\partial}{\longrightarrow}\pi_{0}(H)\to...$  
associated with the principal bundle $H\to G\to G/H$ yields, 
since $G$ is simply connected and $H$ connected, the triviality of the 
fundamental group $\pi_{1}(G/H)=\{1\}$. In other words, 
the manifold $G/H$ is simply connected. One then readily checks that 
the formula 
$$
s_{gH}(g'H)\;=\;g\,\tilde{\sigma}(g^{-1}g')H
$$
defines a ($G$-equivariant) structure $s:G/H\times G/H\to G/H$ of 
symmetric space on the smooth manifold $G/H$. 

\noindent In the above construction, we didn't require the iLa to be 
transvection. This will be used for injectivity of the correspondence.

\noindent Let now $(M,s)$ be a simply connected 
symmetric space and 
denote by $G$ its Lie group of transvections. Choose a base point $o$ in $M$ and 
denote by $H$ its stabilizer in $G$. The natural projection 
$\pi:G\to M:g\mapsto g(o)$ then defines an $H$-principal bundle 
and induces a linear projection
$\pi_{\star e}:\g\to T_{o}(M)$ where $\g$ denotes the Lie algebra of $G$ 
and $T_{o}(M)$ the tangent space to $M$ at point $o$.
Consider now the subgroup $G^{\tilde{\sigma}}\subset G$ 
of fixed points under the involutive automorphism 
$\tilde{\sigma}:G\to G:g\mapsto s_{o}gs_{o}$.
The Lie algebra $\g$ then decomposes as a direct sum of vector space 
$\g=\g^{+}\oplus\CP$ with 
$\tilde\sigma_{\star e}=:\sigma\;=:\;\id_{\g^{+}}\oplus(-\id_{\CP})$.
It then turns out that the Lie algebra $\g^{+}$ of $G^{\tilde{\sigma}}$ 
coincides with the Lie algebra $\h$ of the stabilizer $H$ of $o$.
Indeed, on the one hand, if $Z$ is any element of $\g$, one has 
$\ddto\exp(-t(\sigma(Z)+Z))(o)=Z^{\star}_{o}-Z^{\star}_{o}=0$ 
(where $Z^\star$ is the fundamental vector field associated with $Z\in \mathfrak{g}$)
because\footnote{This follows from the fact that $o$ is an isolated fixed point of 
$s_{o}$} $s_{o_{\star o}}=-\id_{T_{o}(M)}$. Hence the inclusion $\g^{+}\subset\h$. 
On the other hand, for all $h\in H$ and $x\in M$, one has 
$\tilde{\sigma}(h)(x)=s_{o}hs_{o}(x)=s_{o}hs_{o}h^{-1}h(x)=s_{o}s_{h(o)}h(x)
=s_{o}^{2}h(x)=h(x)\,$.
Hence $H\subset G^{\tilde{\sigma}}$. 

\noindent The fact that the obtained iLa $(\g=\h\oplus\CP\,,\,\sigma)$ is 
transvection easily follows from the fact that the transvection group $G$
is the smallest subgroup of automorphisms of $(M,s)$ which is stable under the 
conjugation by $s_{o}$. 

\noindent At last, the canonical isomorphism of symmetric spaces 
$G/H\to M:gH\mapsto g(o)$ yields the announced bijection at the level of 
the isomorphism classes.
\EPf
\begin{rmk}
At this level, one may observe that the Loos connection is the covariant derivative in 
the tangent bundle associated with the canonical connection in the $H$-structure 
$G\to G/H\simeq M$. Indeed, the proof of Theorem \ref{BIJSS} shows that the homogeneous 
space $G/H$ is reductive (the subspace $\CP$ is 
Ad$_H$-invariant in $\g$) in a canonical way. In particular, it is equipped with 
a canonical connection\footnote{sometimes called the Nomizu connection} one form 
(see e.g. \cite{KN1}) $\alpha\in\Omega^{1}(G,\h)$ defined by
$$
\alpha_{g}(\Xi)\;:=\;\pr_{\h}\left(L_{{g^{-1}}_{\star g}}(\Xi)\right)
$$
where $\Xi\in T_{g}(G)$, where $L$ denotes the left-translation action on $G$ 
and where $\pr_{\h}$ denotes the projection from $\g$ onto $\h$ parallel 
to the subspace $\CP$.

\noindent Within that context, the Loos connection $\nabla$ coincides with the 
covariant derivative in the tangent bundle $T(G/H)\simeq G\times_{H}\CP$ associated 
with the canonical connection one form $\alpha$.
\end{rmk}

\noindent Combining Theorems \ref{DERHAM} and \ref{BIJSS},
we immediately get the following geometric decomposition theorem (analogue to the de Rham 
decomposition theorem in Riemannian geometry) which provides a well-defined categorical 
meaning to the Cartesian product in the class of symmetric spaces:
\begin{thm}\label{DERHAMGEO}
Let $(M,s)$ be a simply connected symmetric space. And let 
$$
(M,s)\;=\;\times_{i=1}^{r}(M_{i},s_{i})\;=\;\times_{j=1}^{\overline{r}}
(\overline{M}_{j},\overline{s}_{j})
$$
be two decompositions in indecomposable symmetric spaces. Then, $r=\overline{r}$ and 
there exists a permutation $\tau\in\mbox{Sym}(r)$ and an automorphism  transformation 
$\phi\in\mbox{Aut}(M,s)=\mbox{Aff}(M,\nabla)$
of $(M,s)$ such that for every $i=1,...,r$:
$$
\phi(M_{i})\;=\;\overline{M}_{\tau(i)}\;.
$$
\end{thm}

\subsection{Symplectic symmetric spaces}

\noindent In this section, we consider a class of symplectic manifolds. 
Recall that a symplectic manifold is a pair $(M,\omega)$ where $M$ is a connected 
smooth manifold and where $\omega$ is a non-degenerate closed two-form on $M$. 
A symplectic diffeomorphism, or symplectomorphism, between two symplectic manifolds
$(M,\omega)$ and $(M',\omega')$ is a diffeomorphism $\varphi:M\to M'$ 
such that $\varphi^\star(\omega')=\omega$. In the case $(M,\omega)=(M',\omega')$, 
one speaks about symplectic transformation. The group of symplectic transformations 
is denoted by $Sp(M,\omega)$ or simply $Sp(\omega)$ when no confusion is possible.

\begin{dfn}\label{SSSDEF}\cite{BiThese, Bi2} A {\bf symplectic symmetric space} 
is a triple $(M, \omega,s)$, where $(M,s)$ is a symmetric space in the sense of 
Loos and where $\omega\in\Omega^2(M)$ is a 
non-degenerate two form on $M$ that is invariant under the symmetries: 
for every $x\in M$, one has
$$
s_x^\star\omega\;=\;\omega\;.
$$

\end{dfn}

\begin{dfn} Two symplectic symmetric spaces  $(M,\omega,s)$ 
and $(M',\omega ',s')$ are {\bf isomorphic} if there exists a symplectic 
diffeomorphism $\varphi: (M,\omega) \rightarrow (M',\omega')$ such that 
$\varphi s_x=s'_{\varphi (x)} \varphi$. Such a $\varphi $ is called an {\bf
isomorphism} of $(M,\omega,s)$ onto $(M',\omega',s')$. When 
$(M,\omega,s) = (M',\omega',s')$, one talks about {\bf automorphisms}. 
The group of all automorphisms of the symplectic symmetric space 
$(M,\omega,s)$ is denoted by $Aut(M,\omega,s)$.
\end{dfn}

\noindent Remark \ref{REMPARAL} yields

\begin{prop} On a symplectic symmetric space $(M,\omega,s)$, the 2-form $\omega$
is parallel with respect to the Loos connection:
$$
\nabla\omega\;=\;0\;.
$$
In particular, $\omega$ is closed, hence symplectic.

\noindent Every symplectic symmetric space is therefore a homogeneous 
symplectic space of its transvection group. 
\end{prop}
\noindent In particular, the transvection group $G(M,s)$ of the underlying symmetric 
space (transitively) acts by symplectomorphisms on the 
symplectic manifold $(M,\omega)$. The homogeneous symplectic space underlying a 
symplectic symmetric space is therefore a symplectic cover
of a coadjoint orbit. This aspect is considered in Section \ref{HAM}.

\begin{exs} \label{BiThese,Bi2}
\noindent In dimension two, it turns out that the coadjoint orbit covered by a 
symplectic symmetric space itself admits a structure of symplectic symmetric space. 
We give here, up to isomorphism, the complete list of those two-dimensional coadjoint 
orbits. Four of them are simply connected and two are not.
\begin{enumerate}

\item The flat symplectic plane $(\R^{2}=\{(q,p)\}\,,\,\omega=\lambda\,
{\rm d}p\wedge{\rm d}q)$ ($\lambda\in\R_{0}$) equipped with the Euclidean 
centred symmetries:
$$
s_{x}y\;:=\;2\,x\,-\,y\;.
$$ 
The associated Loos connection if flat, hence metric. The transvection group is the 
translation group $(\R^{2}\,,\,+)$. It is a generic coadjoint orbit of the 
three-dimensional Heisenberg group. It is a Hermitian symmetric space, i.e., 
it admits a compatible complex structure: the one defined by the natural 
identification of $\R^{2}$ with the complex plane $\C$.
\item The two-sphere:
$$
S^{2}_{R}\;:=\;\{x\in\R^{3} \;s.t.\; \lVert x\rVert \;=\;R\}\quad(R>0)
$$
equipped with its standard area form $\omega^{R}$ (induced by the embedding in the Euclidean 
space $\R^{3}$). The symmetry $s_{x}$ at a point $x\in S^{2}_{R}$ is the axial rotation 
of $\pi$ radians around the vector line directed by $\vec{0x}$. 
The Loos connection is the 
Levi-Civita connection w.r.t. the first fundamental form $\beta$ induced by the embedding 
in the Euclidean 3D-space. It is a generic coadjoint orbit of the rotation group $SO(3)$ 
which is its transvection group.
It is a K\"aher symmetric space: the complex structure is the field of endomorphisms 
$J$ of the tangent bundle $T(S^{2}_{R})$ defined by  $\beta(Ju,v):=\omega^{R}(u,v)$. 

\item The hyperbolic plane $\D^{2}_{R}$ of negative curvature $-R^{2}$ $(R>0)$ 
realised as a connected component of an elliptic sphere in Minkowski space-time 
$\R^{(1,2)}$ with scalar product $\eta$ of signature $(-++)$:
$$
\D^{2}_{R}\;:=\;\{x\in\R^{(1,2)}\;|\;\eta(x\,,\,x)\;=\;-R^{2}\}_{0}\;.
$$
The symplectic structure $\omega^{R}$ is the area form induced the natural 
volume form on the ambient Minkowski space. 
The symmetry $s_{x}$ at a point $x\in\D^{2}_{R}$ is the axial Minkowski-rotation 
of $\pi$ radians around the vector line directed by $\vec{0x}$. 
The Loos connection is the Levi-Civita connection w.r.t. the first fundamental 
form ${\beta}$ induced by the embedding in the  Minkowski-space. It is a generic 
coadjoint orbit of the Lorentz group $SO_{0}(1,2)$ which is its transvection group.
It is a K\"aher symmetric space: the complex structure is the field of 
endomorphisms $J$ of the tangent bundle $T(\D^{2}_{R})$ defined by  
$\beta(Ju,v):=\omega^{R}(u,v)$. 

\item The ``massive'' generic co-adjoint orbit $\CO^{m}$ of the Poincar\'e group 
$G:=O_{0}(1,1)\ltimes\R^{(1,1)}$ in dimension $(1,1)$. It is realized as a connected 
component of a generic sphere in $\R^{3}$ equipped with a degenerate scalar product of 
rank two whose isometry group is isomorphic to the Poincar\'e group. 
This degenerate ambient space $\R^{3}$ is isometric to the dual $\g^{\star}$ of 
the Lie algebra $\g$ of the Poincar\'e group equipped with the (dual) Killing metric. 
Within that realization, the orbit $\CO^{m}$ stands as a a cylinder over a branch of 
hyperbola (a mass-shell orbit of mass $m$). The symplectic structure $\omega^{m}$ is then 
given by the Kostant-Kirillov-Souriau 2-form on the coadjoint orbit $\CO^{m}$. Concerning 
the symmetric space structure and associated Loos connection, they can be described as 
follows. Observe first that the connected component $\L$ of the group of affine 
transformations of the real line (the $ax+b$ group) is naturally contained as a Lie 
subgroup in the Poincar\'e group:
$\L\simeq O_{0}(1,1)\ltimes N$ where $N$ is a light ray in the vector subgroup 
$\R^{(1,1)}$. It turns out that the Lie group $\L$ simply transitively acts (by 
restriction of the coadjoint action) on the orbit $\CO^{m}$. 
The affine group $\L$ is exponential in the sense that its exponential map 
realizes a global diffeomorphism between its Lie algebra 
$\fL$ and $\L$. Denoting by $H\in\fL$ an infinitesimal generator of $O_{0}(1,1)$ and by 
$E\in\fL$ an infinitesimal generator of $N$, one gets, for any choice of base point $o$ 
in $\CO^{m}$, a global coordinate system on our orbit:
$$
\phi:\R^{2}\to\CO^{m}\subset\g^{\star}:(a,n)\mapsto\Ad^{\flat}_{\exp(aH)\exp(nE)}(o)\;.
$$
This coordinate system turns out to be a global Darboux chart on $\CO^{m}$:
$$
\phi^{\star}\omega^{m}\;=\;m\,{\rm d}a\wedge{\rm d}n\;.
$$
Note that the map $\phi$ expresses an explicit realization of the orbit 
$(\CO^{m},\omega^{m})$ as the co-tangent bundle (equipped with its canonical Liouville 
symplectic structure) of the above mentioned mass-shell orbit.
The symmetric space structure can now be explicitly described 
within the coordinate system as a ``deformation'' of the flat structure: one has
$$
s_{(a,n)}(a',n')\;=\;\left(2a\,-\,a'\,,\,2\,\cosh(2(a-a'))\,n \,-\,n'\right)\;.
$$
Formula (\ref{CONNLOOS}) easily yields the Loos connection which is not a 
metric connection i.e. it cannot be identified with the Levi-Civita connection 
for any non-degenerate metric on $\CO^{m}$.

\item The generic co-adjoint orbit $\CO^{R}$ of the displacement group of 
the Euclidean plane $G:=SO(2)\ltimes\R^{2}$. 
It is realized as a connected component of a generic sphere in $\R^{3}$ 
equipped with a degenerate scalar product of rank two whose 
group of isometries is isomorphic to $SO(2)\ltimes\R^{2}$. 
This degenerate ambient space $\R^{3}$ is isometric to the dual $\g^{\star}$ 
of the Lie algebra $\g$ of displacement group equipped with the (dual) Killing metric. 
Within that realization, the orbit $\CO^{R}$ stands as a a cylinder over a circle of 
radius $R$. The symplectic structure $\omega^{R}$ is then given by the 
Kostant-Kirillov-Souriau 2-form on the coadjoint orbit $\CO^{R}$. 
It is the Euclidean analogue of the above massive orbit of the Poincar\'e 
group (item 4.). As for the Poincar\'e case, the canonical Loos linear connection on 
$\CO^{R}$ is not a metric connection.

\item Analogue to the hyperbolic plane (item 3.), the generic hyperbolic orbit $\mbox{AdS}_{2}^{R}$ of positive curvature $R^{2}$ of the Lorentz group $SO_{0}(1,2)$:
$$
\mbox{AdS}_{2}^{R}\;:=\;\{x\in\R^{(1,2)}\;|\;\eta(x\,,\,x)\;=\;R^{2}\}\;.
$$
In the same way as for the hyperbolic plane, the anti- de Sitter space 
$\mbox{AdS}_{2}^{R}$ can be realized as a coadjoint orbit of the Lorentz group 
equipped with its first fundamental form (of Lorentz signature in this case) 
induced by the (dual) Killing form on $\mathfrak{so}(1,2)^{\star}$. 
The symmetries are defined in the same way as for the elliptic orbit (item 3.). 
Within this coadjoint orbit realisation, the symplectic structure is the 
canonical one. The transvection group is the Lorentz group $SO_{0}(1,2)$. 
The canonical Loos connection is the Levi-Civita connection associated 
with the Lorentz first fundamental form. It is therefore metric. However, in 
contrast with the hyperbolic plane, the anti-de Sitter space $\mbox{AdS}_{2}^{R}$ 
is para-Hermitian in the sense of Kaneyuki \cite{kaneyuki1985classification} 
(c.f. Definition \ref{PARAKAHLER} below).

\end{enumerate}
It is remarkable that each of those two-dimensional symplectic symmetric spaces is  
canonically associated with some kinematical Lie algebra of space-time dimension two. 
The present article shows that this is not accidental. In fact, all 
the structures described below are somehow already apparent in this 
two-dimensional situation.
\end{exs}

\subsection{Symplectic involutive Lie algebras}
\noindent We started with the notion of involutive Lie algebra 
corresponding to the infinitesimal structure of a (affine) symmetric 
space~\cite{KN1}. 

\noindent Now, we pass to the symplectic context. The tangent version of a 
\emph{symplectic} symmetric space (see \cite{BiThese,BCG2}) is the following one.
\begin{dfn}
A {\bf symplectic iLa} is a triple $(\g\, ,\sigma ,\Omega)$ where $(\g\,,\sigma)$ 
is an iLa and where $\Omega$ is an $\h$-invariant non-degenerate bilinear 
2-form on $\CP$. In particular, the pair $(\CP,\Omega)$ is a \emph{symplectic 
vector space}.

The dimension of ${\cal P}$ defines the {\bf dimension}	of	the triple. 
Two	such triples $(\g_i,\sigma_i,{\Omega}_i)$ $(i=1,2)$ are 
{\bf isomorphic} if there exists a Lie algebra isomorphism 
$\psi :\g_1\rightarrow\g_2$ such that $\psi \circ \sigma_1 = 
\sigma_2 \circ \psi$ and $\psi^*{\Omega}_2={\Omega}_1$. 
\end{dfn}

\noindent Analoguously as for the other types of symmetric spaces, we have

\begin{thm}\cite{BiThese} There is a natural bijective correspondence between the isomorphism
classes of simply connected symplectic symmetric spaces $(M,\omega,s)$ 
and the isomorphism classes of transvection symplectic iLa's $(\g,\sigma,\Omega)$.  
\end{thm}
\noindent The above bijection is easily obtained by  combining  Theorem \ref{BIJSS} with the following basic observation:
\begin{lem}
Let $G$ be a Lie group with Lie algebra $\g$ and $H$ be a closed subgroup (hence Lie and embedded into $G$). Let $\h$ be the Lie algebra of $H$. Then:
\begin{enumerate}
\item[(i)] for every $H$-invariant $(q,p)$-multilinear-tensor $T\in\g/\h^{\otimes q}\otimes((\g/\h)^\star)^{\otimes p}$, the following formula defines a smooth $G$-invariant tensor field $\tilde{T}$ on the smooth manifold $G/H$:
\begin{equation}\label{TENSORINV}
\hspace{-8mm}\left<\tilde{T}_{gH}\,,\,\xi_1\otimes...\otimes\xi_q\otimes v_1\otimes...\otimes v_p\right>\;:=\;
\left<T\,,\,
(g_{\star H})^\star\xi_1\otimes...\otimes(g_{\star H})^\star\xi_q\otimes (g^{-1})_{\star gH}(v_1)\otimes...\otimes (g^{-1})_{\star gH}(v_p)\right>
\end{equation}
where $\xi_i\in T^\star_{gH}(G/H)$ ($i=1,...,q$) and $v_j\in T_{gH}(G/H)$ ($j=1,...,p$) and where 
we canonically identify the tangent space $T_H(G/H)$ to the quotient vector space $\g/\h$.
\item[(ii)] Every smooth $G$-invariant tensor field on the smooth manifold $G/H$ is of the above form.
\end{enumerate}
\end{lem}
\Pf The Lie group $H$ acts on $\g/\h$ by natural projection: 
$h(X+\h)\;:=\;\Ad_h(X)+\h$ ($X\in\g\,,\,h\in H$)
as the vector subspace $\h$ is stable under the adjoint action restricted to $H$. 
One then readily checks that the element $\tilde{T}$ is well-defined as soon as 
it is invariant under this action of $H$ naturally extended to tensors. 
Smoothness is then obvious by the very definition of the smooth structure on $G/H$. 
Item (ii) is immediate from evaluation to the base point $H\in G/H$.
\EPf

\noindent In the symplectic context, one also has a decomposition \`a la de Rham:
\begin{thm}\label{DERHAMS}
Let $(\g,\sigma,\Omega)$ be a transvection siLa. And let 
$$
(\g\,,\,\sigma,\Omega)\;=\;\oplus_{i=1}^r(\g_i\,,\,\sigma_i,\Omega_{i})\;
=\;\oplus_{j=1}^{\overline{r}}
(\overline{\g}_j\,,\,\overline{\sigma}_j,\overline{\Omega}_{j})
$$
be two decompositions into indecomposable siLa's (the direct sum of siLa's 
is defined in the obvious way).
Then, $r\;=\;\overline{r}$ and there exist a permutation $\tau\in\mbox{Sym}(r)$ 
and an automorphism $\varphi$ of $(\g\,,\,\sigma,\Omega)$ such that for 
all $i\in\{1,...,r\}$:
$$
\varphi(\g_i)\;=\;\overline{\g}_{\tau(i)}\;.
$$
\end{thm}

\noindent In the class of symplectic symmetric spaces, one finds those whose 
canonical Loos connection is metric. For instance, K\"ahler (so-called Hermitian) 
symmetric spaces belong to that class. More generally, one has
\begin{dfn}\label{PARAKAHLER}
Let $(\g, \sigma, \Omega)$ be a symplectic triple. One says that it is 
{\bf pseudo-Hermitian} (resp. {\bf para-Hermitian}) if there exists an 
$\h$-commuting symplectic endomorphism $J$ of $\CP$ such that
its square is opposite (resp. equal) to the identity operator on $\CP$. 
Formally, for every $Z\in\h$:
$$
\ad_Z|_\CP\,\circ\,J\;=\;J\,\circ\,\ad_Z|_\CP\quad\mbox{and}\quad 
J^2\;=\;\pm\,\id_\CP\;.
$$
\end{dfn}

\section{Generic kinematical Lie algebras are symplectic involutive}\label{HEART}

\noindent In this section, the ground field $\K$ of our generic kinematical Lie 
algebras is either the reals or the complex field. The first preliminary result 
is the following theorem which proves that to a generic  kinematical Lie algebra 
is naturally associated a symmetric space. 
\begin{rmk}
The proof of Theorem \ref{THM1} below  does not use the requirement on 
the simple $\s$-module $V$ to be absolutely simple.
\end{rmk}
\begin{thm}\label{THM1}
Let $(\g,\s,V)$ be a generic  kinematical Lie algebra. Then, 
the endomorphism $\sigma$ of $\g$ defined by
$$
\sigma\;:=\;\id_{\CZ\oplus\s}\;\oplus\left(-\id_{\CP}\right)
$$
is an involutive automorphism of the Lie algebra $\g$. In other words, 
the pair $(\g\,,\,\sigma)$ is an iLa.
\end{thm}

\Pf
\noindent 
The hypothesis on the isotypical component of $V$ in $\Lambda^{2}(V)$ implies, 
by duality, that there is no nontrivial $\s$-equivariant 
projection from $\Lambda^{2}(V^{\star})$ onto $V^{\star}$. 
Using an $\s$-invariant non-degenerate scalar product on $V$, one concludes that 
there is no nontrivial $\s$-equivariant projection from $\Lambda^{2}(V)$ 
onto $V$. 
This entails that $\s$ cannot non-trivially intersect the \emph{weak}
isotypical component of $V$ in $\g$. Indeed, suppose it does. Then this intersection 
contains a non-trivial simple $\s$-module $V_1$ isomorphic to $V$, in particular 
$V_1$ is an ideal of $\s$ (because it is contained in $\s$). 
So the restriction of the Lie bracket of $\s$ to $V_1$ yields an $\s$-equivariant 
linear map $r$ from $\Lambda^2(V_1)$ to $V_1$. 
But $V_1\cong V$ is $\s$-irreducible, therefore the image of $r$ must be 
isomorphic to either $V$ or trivial.
The hypothesis on isotypical components implies that $V_1$ is an Abelian ideal in $\s$. 
As this Abelian ideal is isomorphic to $V$ as an $\s$-module this would contradict the 
hypothesis of faithfulness  of $V$ as a module of $\s$.

\vspace{1mm}

\noindent In the same way, considering a decomposition of $\CP$ into isomorphic 
simple $\s$-modules (each of them isomorphic to $V$):
$$
\CP\;=\;\CP_{0}\;\oplus\;\CP_{1}\;,
$$
for all $i,j=0,1$, the $\s$-module $[\CP_{i}\,,\,\CP_{j}]$ cannot intersect 
the isotypical component of $V$ in $\g$. 

{\color{blue}

\noindent Indeed, let us first choose two isomorphisms of $\s$-modules:
$$
\phi_k:V\to\CP_k:v\mapsto v_k\quad(k=0,1)\;.
$$
Consider the decomposition into $\s$-submodules:
$$
\Lambda^2(\CP)\;=\;\Lambda^2(\CP_0)\;\oplus\;\Lambda^2(\CP_1)\;\oplus\;\CP_0\wedge\CP_1
$$
where $\CP_0\wedge\CP_1$ is the linear subspace of $\Lambda^2(\CP)$ generated by 
the elements of the form $p_0\wedge p_1'$ with $p_0\in\CP_0$ and $p_1'\in\CP_1$.

\noindent Denote by
$$
C:\Lambda^2(\CP)\to\g
$$
the homomorphism of $\s$-modules given by 
$$
C(p\wedge p')\;:=\;[p,p']\;.
$$

\noindent Now, given $i$ and $j$, define the homomorphism of $\s$-modules
$$
\Delta_{ij}:\Lambda^2(V)\to\Lambda^2(\CP):v\wedge w\mapsto v_i\wedge w_j\;.
$$

\noindent Denote by 
$$
\pr_i:\g\to\CP_i
$$
the natural projection parallel to $\h\oplus\CP_{\overline{i}}$ with $\overline{i}\neq i$. And finally consider the homomorphisms of $\s$-modules
$$
\Delta_{ijk}:\Lambda^2(V)\to  V
$$
defined by
$$
\Delta_{ijk}\;:=\;\phi_k^{-1}\,\circ\,\pr_k\,\circ\,C\,\circ\,\Delta_{ij}\;.
$$
The hypothesis on the isotypical component of $V$ in $\Lambda^2(V)$ implies that 
these maps $\Delta_{ijk}$ are all identically zero. Therefore, observing that 
$$
\phi_k^{-1}\left(\pr_k[v_i,w_j]\right)\;=\;\Delta_{ijk}(v\wedge w)\;,
$$
one obtains}
the inclusion
$$
[\CP\,,\,\CP]\;\subset\;\h\;:=\;\CZ\,\oplus\,\s\;.
$$

\noindent Let us now analyse the space $[\CZ,\CP]$. 
If, for some $j=0,1$, the $\s$-submodule $[\CZ,\CP_{j}]$ is nontrivial, 
it must be isomorphic to $V$ as an $\s$-module. What has been observed in the first 
preceding paragraph combined with a dimensional argument, $[\CZ,\CP_{j}]$  therefore 
lives in the isotypical component of $V$ in $\g$, that is, it lives in $\CP$.
Hence the structure of involutive Lie algebra.\EPf

\noindent We end this section by introducing a map 
$\Omega:\CP\times\CP\to \K\,$, 
with no restriction on the ground field ($\K=\R$ or $\C$), 
which will be crucial in the rest of this article, together with 
the following
\begin{lem} The $\K$-bilinear map
\[
\Omega:\CP\times\CP\to\h/\s\;\simeq\;\K:(X,Y)\mapsto [X,Y]\;+\;\s\;
\]
is $\h$-equivariant. 
\end{lem}
\Pf Since $\CZ$ is central in $\h$, for all $H\in\h$ and $X,Y\in\CP$, the iLa structure 
implies
$$
\Omega([H,X]\,,\,Y)\;+\;
\Omega(X\,,\,[H,Y])\;=\;\left[[H,X]\,,\,Y\right]\;+\;\left[X\,,\,[H,Y]\right]\;
+\;\s\;.
$$
Then, by using the Jacobi identity and that fact that 
$\CZ$ is central in $\h$, one finds
$\left[[H,X]\,,\,Y\right]\;+\;\left[X\,,\,[H,Y]\right]\equiv$  
$-\left[H,[X,Y]\right]\in [\h,\h]\subset \s$. By definition of $\Omega\,$, 
we therefore obtained that $\Omega([H,X]\,,\,Y)\;+\;\Omega(X\,,\,[H,Y])\;= 0$, 
which proves the $\h$-invariance of $\Omega$.
\EPf

\noindent The iLa structure being established, we can now pass to the additional 
symplectic structure. We will first consider the case of a complex generic 
kinematical Lie algebra.
 
\subsection{The complex case}
\noindent In this paragraph, we consider the case where the base field $\K=\C$ 
is the one of complex numbers. It means that the Lie algebras and modules are 
all over the complex numbers. 
\begin{rmk}
The proofs of the results in this paragraph  do not use the requirement 
on the simple $\s$-module $V$ to be absolutely simple.
\end{rmk}

\begin{thm}\label{COMPLEXSILA} Let $(\g,\s,V)$ be a complex generic kinematical 
Lie algebra. Then

\begin{enumerate}

\item[(i)] if the  (complex bilinear) map
$$
\Omega:\CP\times\CP\to\h/\s:(X,Y)\mapsto[X,Y]\,+\,\s
$$
has a trivial radical, then the triple $(\g,\sigma,\Omega)$ is a (complex) siLa 
(where we (non-canonically) identify  $\h/\s$ to the field of complex numbers $\C$).

\item[(ii)] If the above map $\Omega$
admits a nontrivial radical, then the symmetric space associated with the iLa 
$(\g,\sigma)$ is decomposable.

\item[(iii)] When decomposable, the symmetric space is flat. In particular, 
since it is even-dimensional and flat, it admits a structure of symplectic 
symmetric space.

\end{enumerate}
\end{thm}

\Pf 
\noindent The radical of $\Omega$
$$
D\;:=\;\mbox{rad}(\Omega)\;:=\;\{X\in\CP\;|\;\Omega(X,Y)\;=\;0\;,\,\forall\, Y\in\CP\}
$$
is an $\h$-submodule of $\CP$ -- for $X\in \mbox{rad}(\Omega)$ 
and $H\in \h$, 
one has $\Omega([H,X],Y)=-\Omega(X,[H,Y])$ which vanishes because 
$[H,Y]\in \mathcal P$ and $X\in \mbox{rad}(\Omega)$, 
hence $[H,X]\in \mbox{rad}(\Omega)$ whenever $X\in \mbox{rad}(\Omega)$ --  
such that
\begin{equation}\label{DPS}
[D\,,\,\CP]\;\subset\;\s\;,
\end{equation}
by definition of $\Omega$. As a result, the quotient space $\CP/D$ is an 
$\h$-module.

\noindent As an $\s$-submodule, $D$ must either be reduced to zero, equal to $\CP$ 
or 
isomorphic to $V$.
\begin{enumerate}
\item Suppose $D\cong V$. 
In that case, the $\h$-module $\CP/D$ is isomorphic to $V$. 
Let us consider the $\s$-commuting action of $\CZ\subset \h$ 
on the $\h$-module $\CP/D$. Since the latter is a simple $\s$-module 
isomorphic to $V$, Schur's lemma implies that the latter action must 
be proportional to the identity. 
On the other hand, since the $\h$-module $\CP/D$ is symplectic
--- for any $X\in\CP/D$, the condition 
$\Omega(X,Y)=0$ $\forall$ $Y\in \CP/D$ implies that $X=0$, so that 
$\Omega\vert_{\CP/D \times \CP/D}$ is non-degenerate; moreover, the 
the $\h$-invariance of $\Omega$ is carried to $\CP/D \times \CP/D$ ---$\,$,
we conclude that the proportionality constant must be zero, since 
the identity map cannot be symplectic.
In particular, for every $\s$-module $W$ supplementary to $D$ in $\CP$, 
one must have
$$
[\CZ\,,\,W]\;\subset\;D\;.
$$
Still in the situation where $D$ is isomorphic to $V$, let $d\in D$. 
Then, since $[W,D]$ lives in $\s$ -- recall \eqref{DPS}--, the identity
\begin{eqnarray*}
&&
0\;=\;[w_{1}\,,\,[w_{2}\,,\,d]]\;+\;[d\,,\,[w_{1}\,,\,w_{2}]]\;
+\;[w_{2}\,,\,[d\,,\,w_{1}]]
\end{eqnarray*}
implies that the middle term lives in $D$ while both the other ones live in $W$. 
Therefore, we have the relation:
\begin{equation}\label{WWD}
[[W\,,\,W]\,,\,D]\;=\;\{0\}\,.
\end{equation}
On the other hand, since $D$ is an $\h$ (hence $\s$) module 
and because $\CZ$ commutes with $\s$, Schur's lemma implies that the 
action of $\CZ$ on $D\simeq V$ is scalar. 
Considering the components 
$$
[W\,,\,W]\;\subset\;[W\,,\,W]_{\CZ}\;\oplus\;[W\,,\,W]_{\s}
$$
according to the above decomposition of $\h$, we conclude, from 
the inclusion (\ref{WWD}), that the action of $[W\,,\,W]_{\s}$ on $D$ must 
be scalar too, in the sense that there exists a one form $\alpha$ on $[W,W]_\s$ 
such that for every $T\in[W,W]_\s$, one has
$$
\ad_T|_D\;=\;\alpha(T)\,\id_D\;.
$$

\noindent Now, as $\s$-modules, $D$, $W$ and $\CP/D$ are isomorphic. 
And on the last one, there exists, by construction, an $\s$-invariant 
symplectic structure. So, such an $\s$-invariant symplectic 
structure exists on the three of them (possibly different from the restriction of 
$\Omega$). Therefore,
for every $T$, one must have $\alpha(T)=0$, i.e., $\alpha$ must be identically zero. 

\noindent Now, by the fact that $D\cong V$ is acted upon faithfully by $\s$, 
the relation $\ad_{[W,W]_\s}\vert_D=0$ implies that $[W,W]_\s = \{0\}$, hence 
$[W,W]\subset {\cal Z}$. Since $\CZ$ is one-dimensional, one has either 
$[W,W] = \CZ$ or $[W,W] = \{0\}$. But since $W$ is symplectic w.r.t. $\Omega$,
we get the equality:
$$
[W\,,\,W]\;=\;\CZ
$$
implying from \eqref{WWD} that
$$
[\CZ,D]\;=\;\{0\}\;.
$$
We recall that we had obtained $[\CZ,W]\subset D\,$ earlier.

\noindent We will now prove, still under the hypothesis that $D$ is isomorphic to $V$, 
that $[\CP\,,\,\CP]\subset\CZ$.
In order to do so, we first remark that $D$ is Abelian. 
Indeed, Jacobi tells us that
$$
[[D\,,\,D]\,,\,W]\;\subset\;[[W,D]\,,\,D]\;\subset\;D\;.
$$
Since $[D,D]\subset\s$, one also has  $[[D\,,\,D]\,,\,W]\subset W$. Hence 
$$
[[D\,,\,D]\,,\,W]\;\subset\;D\cap W\;=\;\{0\}\;.
$$
From the faithfulness of $W$ as a $\s$-module, one must have $[D,D]=\{0\}$.

\noindent To further proceed in proving the inclusion $[\CP\,,\,\CP]\subset\CZ$, 
we now consider the following contraction of our iLa 
by considering the new bracket $[.\,,\,.]_{0}$ on $\g_{0}\;:=\;\g$ defined 
by (with obvious notations):
$$
\begin{array}{c}
\left[W\,,\,W\right]_{0}\;:=\;\{0\}\\
\left[\CZ\,,\,\g_0\right]_{0}\;:=\;\{0\}\\
\left[\g_0\,,\,D\right]_{0}\;:=\;\left[\g\,,\,D\right]\\
\left[\s\,,\,\g_0\right]_{0}\;:=\;\left[\s\,,\,\g\right]\;.
\end{array}
$$
The space $\g_{0}\;=\;\h\;\oplus\;\CP$ again underlies an iLa. 
Indeed, consider $d\in D$ and $w_{1},w_{2}\in W$. Then, using (\ref{DPS}) and 
(\ref{WWD}):
\begin{eqnarray*}
&&
0\;=\;[w_{1}\,,\,[w_{2}\,,\,d]]\;+\;[d\,,\,[w_{1}\,,\,w_{2}]]\;+\;[w_{2}\,,\,
[d\,,\,w_{1}]]\;=\;[w_{1}\,,\,[w_{2}\,,\,d]]\;+\;[w_{2}\,,\,[d\,,\,w_{1}]]\\
&=&
[w_{1}\,,\,[w_{2}\,,\,d]_{0}]\;+\;[w_{2}\,,\,[d\,,\,w_{1}]_{0}]
\end{eqnarray*}
which equals
$$
[w_{1}\,,\,[w_{2}\,,\,d]_{0}]_{0}\;+\;[w_{2}\,,\,[d\,,\,w_{1}]_{0}]_{0}\;
=\;[w_{1}\,,\,[w_{2}\,,\,d]_{0}]_{0}\;+\;[w_{2}\,,\,[d\,,\,w_{1}]_{0}]_{0}
\;+\;[d\,,\,[w_{1}\,,\,w_{2}]_{0}]_{0}
$$
since $[W,D]$ lives in $\s$.
The remaining Jacobi identities are immediate.

\noindent Now, since $\k_{0}:=[\CP\,,\,\CP]_{0}\subset\s$ and $\s$ 
acts faithfully on $\CP$, the space
$\k_{0}\;\oplus\;\CP\subset\g_{0}$ is the transvection iLa of a 
(pseudo-)Riemannian symmetric space.
The decomposition $\CP\;=\;W\oplus D$ is stable under the action of the holonomy 
Lie algebra $\k_{0}$, hence as a consequence of the de Rham-Wu decomposition 
theorem (\cite{wu1963rham}), this symmetric space is reducible, entailing
$$
[D,W]\;=\;[D,W]_{0}\;=\;\k_{0}\;=\;\{0\}\;.
$$
This means that the pseudo-Riemannian symmetric space at hand is flat.

We now resume to our analysis of the kinematical algebra $\g$ 
and cease considering $\g_0$.
We choose a generator $Z_0$ of $\CZ$ and consider the $\s$-intertwiner
$$
\rho\;:=\;\ad_{Z_0}|_W\;:W\to D\;.
$$
\noindent Consider the symplectic vector space $(W,\Omega^W)$ defined by
$$
[w,w']\;:=\;\Omega^W(w,w')\,Z_0\;.
$$
\noindent Then, the Jacobi identity tells us that
$$
0\;=\;\oint_{1,2,3}[w_1\,,\,[w_2\,,\,w_3]]\;
=\;-\,\oint_{1,2,3}\Omega^W(w_2,w_3)\,\rho(w_1)
$$
which non-degeneracy of $\Omega^W$ forces
$$
\rho\;\equiv\;0\;,
$$
which shows that $[\CZ,W]=\{0\}\,$.
Above, we had found that $[\CZ,D]=\{0\}\,$, hence we have that 
$[\CZ,\CP]=\{0\}\,$.
\noindent The Lie algebra $\g$, in the case 
$\mbox{rad}(\Omega)\;\cong\;V$, is therefore the semi-direct product of
$\s$ acting on a decomposable flat iLa  whose one direct factor is a 
Heisenberg Lie algebra $(\CZ\;\oplus W)$:
$$
\g\;=\;\s\ltimes\left(D\oplus(\CZ\;\oplus W)\right)\;.
$$

\item
In case the radical $\rad(\Omega)$ is the entire space $\CP$, i.e. 
$\fK:=[\CP,\CP]\subset\s$, the space $\fK\oplus\CP$ is the transvection iLa of a 
holonomy reducible (pseudo-)Riemannian symmetric space \cite{deRhamWu}. 
In particular, the holonomy Lie algebra $\fK$ splits into a direct sum of ideals:
$$
\fK\;=\;\fK_{0}\;\oplus\;\fK_{1}
$$
that are such that, say, $[\fK_{0}\,,\,\CP_{1}]\;=\;\{0\}$, contradicting the 
faithfulness of $\CP_{1}\,\simeq\,V$ as soon as $\fK$ is nontrivial. 
Therefore, since $\s$ acts faithfully on $\CP$, 
we conclude that $\fK=[\CP,\CP]=\{0\}$. We therefore have 
$$
\g\;=\;\s\ltimes\left(\CP\oplus\CZ\right)\;.
$$
The only constraint we have on $\CZ$ is that it is a line in the centralizer 
of $\s$, in the endomorphisms of $\CP$.

All the possibilities for the action of $\CZ$ on $\CP$ are possible. 

\item It remains to analyse the case where $(\CP,\Omega)$ is symplectic, i.e., 
the iLa underlies a symplectic symmetric space
modelled on $(\CP,\Omega)$. We assume here that the associated symmetric space is 
decomposable.

\vspace{2mm}

\noindent Let us first consider the case where $\CZ$ is central in $\g$. 
In this case, the holonomy $\fK$ of the transvection Lie algebra is a subalgebra of 
$\s$ (because $\s$ faithfully acts on $\CP$). 
On the other hand, for decomposability, there exists a $\fK$-stable decomposition 
$\CP\;=\;\CP_0\oplus\CP_1$ with commuting factors $[\CP_0\,,\,\CP_1]=\{0\}$. 
Note that $[\CP_0\,,\,\CP_0]\subset\fK$ trivially acts on $\CP_1$, contradicting the 
faithfulness of the action of $\s$ on $V\simeq\CP_1$ if non-trivial. And similarly for 
$[\CP_1\,,\,\CP_1]$. Both factors are therefore flat, i.e. 
$[\CP_0\,,\,\CP_0]=[\CP_1\,,\,\CP_1]=\{0\}$. Hence $\fK=\{0\}$.

\vspace{2mm}

\noindent Now, let us consider the case where $\CZ$ non-trivially acts on $\CP$. 
In that case, the holonomy is of the form:
$$
\fK\;=\;\CZ\oplus\fK_\s
$$
where $\fK_\s\subset\s$. Under the hypothesis of decomposability, the same argument 
as above implies
$\fK_\s=\{0\}$. Which yields, say, $[\CP_0\,,\,\CP_0]=\CZ$ hence, necessarily, 
$[\CP_1\,,\,\CP_1]=\{0\}$,
in contradiction with the fact that $\Omega$ is symplectic.
\EPf
\end{enumerate}

\noindent We end this section by observing the following dichotomy concerning the 
action of the (complex) line $\CZ$ on $\CP$:

\begin{prop}\label{PROPNORS}
The action of the (complex) line $\CZ$ on $\CP$ is either nilpotent or semisimple.
Moreover, when nilpotent, the action of $\CZ$ on $\CP$ squares to zero.
\end{prop}
\Pf Let $Z_0$ be a generator of the complex line $\CZ$ and denote by $A$ its 
(complex linear) action on $\CP$:
$$
A\;:=\;\ad_{Z_0}|_{\CP}\;:\CP\to\CP\;.
$$
The endomorphism $A$ uniquely decomposes (Jordan decomposition) into a sum of 
two commuting endomorphisms:
\begin{equation}\label{JORDAN}
A\;=:\;S\;+\;N
\end{equation}
where $S$ is semisimple and $N$ nilpotent. Remind that, as $\Omega$ is $\h$-invariant, 
the endomorphism $A$ is symplectic: $A\in\sp(\CP,\Omega)$.

\noindent Assume the kernel $\CN$ of $N$ to be a proper subspace of $\CP$ 
(i.e. $\{0\}\neq\CN\neq\CP$). The operator $S|_\CN$ acts as $A$ on $\CN$ and therefore 
commutes with $\s$. Hence it must be scalar:
$$
S|_\CN\;=:\;\lambda\,\id_\CN\quad(\lambda\in\C).
$$
Now, if $S$ admits another eigenvalue, say $\mu\neq\lambda$, then denoting by 
$\CP_\mu$ the corresponding eigenspace in $\CP$, one has 
$$
N(\CP_\mu)\;\subset\;\CP_\mu\;.
$$
But $N$ cannot admit a kernel in $\CP_\mu$ because that one would be contained in 
$\CN\subset\CP_\lambda$. The operator 
$N$ being nilpotent, this  yields a contradiction. Hence, when $\CN$ is proper in 
$\CP$, the semisimple part $S$ of $A$ must acts as a scalar on the entire $\CP$. 
Now, since $\sp(\CP,\Omega)$ is semisimple, it contains both nilpotent and semisimple 
components of any of its elements. Therefore, $S$ must be identically zero (because 
scalar).

\noindent At last, observe that the kernel $\CN$ of $N=\ad_{Z_0}|_\CP$ is invariant 
under the action of $\h$. If symplectic, the argument in the proof of Proposition 1.2 
page 249 of \cite{Bi98bis} implies that it commutes with its orthogonal subspace in 
$\CP$. Hence Theorem \ref{THM1} item (ii) implies flatness.

\noindent The image $\CI:=N(\CP)$ of $N=\ad_{Z_0}|_\CP$ is $\s$-invariant. 
Indeed, if $N(X)$ ($X\in\CP$) is an element of $\CI$, we have for every $T\in\s$: 
$[T,N(X)]=[T,[Z_0,X]]=[Z_0,[T,X]]=N([T,X])\in\CI$. The subspace $\CI$ therefore 
cannot properly intersect $\CN$, because this intersection, whose intersection is 
strictly smaller than the one of $\CN$, would then belong to another isotypic component 
than the one of $V$ in $\CP$. This implies
$\CI=\CN$ i.e. $N^2=0$.
\EPf

\begin{prop}\label{IRREDCOMPLEX}
\noindent Let us assume the complex symplectic iLa  $(\g,\sigma,\Omega)$ to be 
indecomposable and non-flat. Assume the endomorphism $A$ to be semisimple.
Then:
\begin{enumerate}
\item[(i)] the space $\CP$ decomposes into a direct sum of two $\h$-invariant 
transverse Lagrangian subspaces in duality:
$$
\CP\;=\;L\;\oplus\;\overline{L}\;.
$$
\item[(ii)] Both $L$ and $\overline{L}$ are Abelian subalgebras of $\g$.
\item[(iii)] Both $L$ and $\overline{L}$ are simple $\s$-modules.
\item[(iv)] The endomorphism $A$ acts on $L$ and $\overline{L}$ with opposite 
eigenvalues.
\end{enumerate}
\end{prop}

\Pf  Since, in this case, $A$ is semisimple, the space $\CP$ decomposes into $\CZ$-root spaces:
$$
\CP\;:=\;\oplus_{\alpha\in\Phi}\CP_\alpha
$$
where $\Phi\subset\CZ^\star$. Each of these root spaces being $\h$-invariant, the cardinality of $\Phi$ is at most two.
In the case it is equal to two, the argument in the proof of Proposition 1.2 page 249 of \cite{Bi98bis}  yields the assertion.

\noindent At last, the endomorphism $A$ cannot be scalar as it would not be symplectic in that case.
For the same reason, in view of the matrix form of  $\Omega$ in the decomposition (i), one gets item (iv).

\noindent Note also that the map $\overline{L}\to L^\star$ canonically defined by the bilinear two-form $\Omega$ admits a non-trivial kernel only if the radical of $\Omega$ is non-trivial for both subspaces are isotropic. It is not the case in the present situation, hence the Lagrangian subspaces $L$ and $\overline{L}$ are in duality.
\EPf

\noindent {\bf Warning:} From now on, all the symmetric spaces considered below will be 
\underline{indecomposable and non-flat}.

\subsection{The real case}

\noindent Let us now come back to the case where the ground field $\K=\R$ is the reals. 
In this real case, \emph{we restore the hypothesis on the real $\s$-module $V$ to be 
absolutely simple}. 
We start by observing 
\begin{prop} Let $(\g,\s,V)$ be a real kinematical Lie algebras with indecomposable and non-flat associated iLa $(\g,\sigma)$. Then, the iLa 
 $(\g,\sigma)$ canonically underlies a symplectic iLa $(\g,\sigma,\Omega)$.
\end{prop}
\Pf 
We consider the element $\Omega:\CP\times\CP\to\h/\s:(X,Y)\to[X,Y]+\s$.
As $V$ is an absolutely simple real $\s$-module, we have seen that the complexified iLa 
$\g^\C=\h^\C\oplus\CP^\C$ associated with the complex generic kinematical 
$(\g^\C,\s^\C,V^\C)$ Lie algebra becomes symplectic when equipped with the complexified 
element $\Omega^\C$. This implies that the  radical of the real element $\Omega$ must be 
trivial.
\EPf

\subsubsection{The action of $\CZ$ on $\CP$ is semisimple}

\noindent According to the dichotomy resulting from Proposition \ref{PROPNORS}, we first investigate the case where the action of $\CZ$ on $\CP$ to be (possibly complex) semisimple.

\noindent As earlier, we consider a generator $Z_0$ of the $\K$-line $\CZ$ and denote by $A$ its action on $\CP$:
$$
A\;:=\;\ad_{Z_0}|_{\CP}\;:\CP\to\CP\;.
$$

\vspace{2mm}

\noindent The complex linear extension $A^\C$ of $A$ on  the complexified space $\CP^\C$ is semisimple and one has the decomposition  into $\CZ^\C$-root spaces:
$$
\CP^\C\;:=\;\CP_\alpha\;\oplus\;\CP_{-\alpha}
$$
where $\alpha\in(\CZ^\C)^\star$. 

\noindent Denoting by $X\mapsto\overline{X}$ the conjugation of $\g^\C$ 
associated with the real form $\g$, one observes that for every $Z$ in $\CZ\subset\CZ^\C$:
$$
\overline{[Z\,,\,X_\alpha]}\;=\;\overline{\alpha(Z)}\,\overline{X_\alpha}\;=\;[\overline{Z}\,,\,\overline{X_\alpha}]\;=\;[Z\,,\,\overline{X_\alpha}]\;.
$$
Hence $\overline{X_\alpha}$ is an $A$-eigenvector, implying that $\lambda:=\alpha(Z_0)$ is either real or purely imaginary.

\begin{lem} 
In the case $\lambda$ is real, $A$ is real-semisimple and $\CP$ admits an $\h$-invariant para-complex structure that is realized by the action of an element of $\CZ$.
\end{lem}
\Pf We set $X_\alpha\;=:\;u_\alpha\;+\;i\,v_\alpha$ with $u_\alpha$ and $v_\alpha$ in $\CP$. 
We then observe:
$$
A(X_\alpha)\;=\;\lambda(X_\alpha)\;=\;\lambda u_\alpha\;+\;i\lambda\,v_\alpha\;=\;A(u_\alpha)\;+\; i\,A(v_\alpha)\;,
$$
hence 
$$
A(u_\alpha)\;=\;u_\alpha\;.
$$
And similarly for $\CP_{-\alpha}$. \EPf

\noindent Following the same lines, we get
\begin{lem} 
In the case $\lambda$ is purely imaginary, $\CP$ admits an $\h$-invariant complex structure that is realized by the action of an element of $\CZ$.
\end{lem}

\noindent As an immediate corollary, we have

\begin{prop}
Let $(\g,\s,V)$ be a real generic kinematical Lie algebra with indecomposable 
non-flat associated iLa $(\g,\sigma,\Omega)$. Suppose the action of $\CZ$ on $\CP$ 
to be semisimple. Then:
\begin{enumerate}
\item[(i)] if the action of $\CZ$ is not real semisimple then the symmetric space 
admits a symplectic-compatible complex structure:
the siLa is pseudo-Hermitian.
\item[(ii)] If the action of $\CZ$ is real semisimple then the symmetric space 
admits a symplectic-compatible para-complex structure;
the siLa is para-Hermitian.
\end{enumerate}
\end{prop}
\noindent The above result classifies the generic kinematical Lie algebras 
where $\s$ is semisimple. Indeed, one has
\begin{thm}\label{SSSZSS}
Let $(\g,\s,V)$ be a generic kinematical Lie algebra over $\K=\R$ or $\C$ such 
that $\s$ is semisimple and the action of $\CZ$ on $\CP$ is semisimple. Assume the 
associated siLa to be indecomposable and non-flat. Then 
\begin{enumerate}
\item[(i)] $\g$ is simple. 
\item[(ii)] The associated siLa $(\g,\sigma,\Omega)$ is the transvection siLa of a 
simple symplectic symmetric space.
\item[(iii)] The Lie algebra $\s$ is compact if and only if the symplectic symmetric 
space is Hermitian.
\item[(iv)] The Lie algebra $\s$ is complex if and only if the symplectic symmetric 
space is hyper-K\"ahler.
\item[(v)] The Lie algebra $\s$ is not complex nor compact if and only if the 
symplectic symmetric space is causal of Cayley type
(para-Hermitian).
\end{enumerate}
\noindent Reciprocally, every simple symplectic iLa is the associated siLa with a generalized a generic kinematical Lie algebra such that $\s$ is semisimple and the action of $\CZ$ on $\CP$ is semisimple.
\end{thm}
\Pf The fact that the action of $\h$ is reductive on $\CP$ forces $\g$ to be a reductive Lie algebra. By indecomposability of the symmetric space, it must be simple \cite{BiThese, Bi98bis}. Items (ii)-(v) then follows from \cite{BiThese, Bi98bis}.
Regarding the last assertion is structural: it immediately follows from the fact that the simple siLa are the simple iLa whose Levi factor of the holonomy act reducibly on $\CP$.
\EPf

\subsubsection{The Levi condition}

\noindent Motivated by the classical notion of kinematical Lie algebra, i.e., $\s\simeq\mathfrak{so}(D)$, we observe
\begin{prop}\label{COMPACTINLEVISS} Let $(\g,\sigma)$ be an iLa and 
$\s$ be a compact and semisimple Lie subalgebra of $\g$. Then it is contained in 
a $\sigma$-stable Levi factor of $\g$.
\end{prop}

\noindent The proof follows from a few observations.
First, we have

\begin{prop}
Let $K$ be a connected compact real Lie group. Let $K=TK_0$ be a decomposition of $K$ into a Lie group direct product of its centre $T$ and a compact semisimple Lie group $K_0$.

\noindent  Let $S$ be a semisimple Lie subgroup of $K$ with no centre. 
Then $S$ is contained in $K_0$, which is in particular unique.
\end{prop}

\Pf
It is sufficient to assume $S$ simple (apply the argument below to every of its simple components). Let 
$
\tau:S\to T
$ be
the map defined by the global decomposition $K=TK_0$. The decomposition $K=TK_0$ is a direct product of groups hence the map $\tau$ is a Lie group homomorphism.
Therefore its kernel is a normal subgroup of $S$ and one concludes by simplicity.
\EPf
\noindent Therefore one has the
\begin{cor}\label{COMPACTINLEVI}
Every compact semisimple Lie subalgebra $\s$ of a real Lie algebra $\CG$ is contained in a Levi factor of $\CG$.
\end{cor}
\Pf
Let $\fK_0$ be a maximal compact Lie subalgebra of a Levi factor $\fL$ of $\CG$. And let $\fK$ be a maximal compact Lie subalgebra of $\CG$ containing $\fK_0$. Let $\CR$ be the solvable radical of $\CG$ and consider the projection $\pi_\fL$ from $\CG$ onto $\fL$ parallel to $\CR$. The map
$\pi_\fL$ is a Lie algebra homomorphism. Therefore $\pi_\fL(\fK)$ is a compact subalgebra in $\fL$, therefore contained in $\fK_0$. Since $\fK$ contains $\fK_0$, one has the equality 
$\fK_0=\pi_\fL(\fK)$ and the splitting
$$
\fK\;=\;\fK_0\oplus(\fK\cap\CR)\;.
$$
Indeed, for every $k\in\CK$, one has $k=\pi_\fL(k)\oplus\pi_\CR(k)$. Since $\pi_\fL(k)$ belongs to $\CK$ as we proved above, we have $\pi_\CR(k)=k-\pi_\fL(k)$ belongs to $\CK$ and $\CR$. Hence $\pi_\CR(\CK)=\CK\cap\CR$.

\noindent The Levi decomposition of $\fK$ is
$$
\fK\;=\;\fK_0'\oplus\fZ_0\oplus(\fK\cap\CR)
$$
where $\fZ_0$ is the center of $\fK_0$. Hence $\fK\cap\CR$, being a solvable ideal of a compact Lie algebra $\fK$ must be central in $\fK$. And by the preceding proposition, $\s$ is, up to conjugation\footnote{This was first proven by Iwasawa and Cartan, but we refer to \cite{borel1950sous} for a self-contained presentation.}, contained in $\fK_0'$.
\EPf
\noindent Proposition \ref{COMPACTINLEVISS} now follows from the fact that every iLa $(\g,\sigma)$ admits a $\sigma$-stable Levi factor $\CL$, necessarily conjugated to $\fL$.

\noindent This observation leads us to consider a special class of generic kinematical Lie algebras:
\begin{dfn}\label{LEVIC}
A generic kinematical Lie algebra $\g$ satisfies the {\bf Levi condition} is $\s$ is contained in a $\sigma$-stable Levi factor $\CL$ of $\g$.
\end{dfn}

\subsubsection{The $\CZ$-nilpotent case : Poincar\'e spaces}\label{POINCARESECT}

\noindent  We now
investigate the case where the action of $\CZ$ on $\CP$ is nilpotent. Remind that in this case, one has $A^2=0$
(c.f. the proof of Proposition \ref{PROPNORS}).

\begin{dfn}\label{POINCARE}
A generic kinematical (real or complex) Lie algebra is called of the {\bf Poincar\'e type} if 
\begin{enumerate}
\item it satisfies the Levi condition \ref{LEVIC}.
\item The action of $\CZ$ on $\CP$ is nilpotent.
\item The underlying symmetric space is not solvable (i.e. its transvection iLa contains a nontrivial Levi factor).
\item The underlying symmetric space is indecomposable.
\end{enumerate}
\end{dfn}

\begin{prop} Let $(\g,\sigma,\Omega)$ be a siLa underlying a Poincar\'e generic kinematical Lie 
algebra over $\K=\R$ or $\C$. Then:
\begin{enumerate}
\item[(i)] the radical $\CR$ is Abelian.
\item[(ii)] The subspaces $\CP_{\CR}\;:=\;\CP\cap\CR$ and $\CP_{\CL}:=\CP\cap\CL$ of $\CP$ are Lagrangians in duality.
\item[(iii)] The subspaces $\CP_{\CR}$ and $\CP_{\CL}$ of $\CP$ are $\s$-modules each of them isomorphic to $V$.

\item[(iv)] One has the equality:
$$
\left[\CP_{\CR}\,,\,\CP_{\CL}\right]\;=\;\CZ\;.
$$
\item[(v)] The action of $\CZ$ establishes an isomorphism of $\s$-modules from $\CP_{\CL}$ to $\CP_{\CR}$.

\end{enumerate}
\end{prop}

\Pf Since the action of $\CZ$ is nilpotent,  $\g$ is not semisimple. 
The subspace $\CP_{\CR}$ is non-trivial and different from $\CP$ (otherwise 
the symmetric space would be solvable). 
Therefore it is an $\s$-submodule of $\CP$ isomorphic to $V$. 

\noindent The subspace $[\CP,\CP_{\CR}]$ is non-trivial (as $\Omega$ is non-degenerate) and contained in $\h\cap\CR$. Hence it cannot intersect $\s$ and must therefore be equal to $\CZ$. The line $\CZ$ must therefore live in $\CR$, entailing $\CR\cap\h=\CZ$ and $\h\cap\CL=\s$. Moreover, $\CP_{\CR}$ cannot be symplectic, because of indecomposability (Proposition 4.1  page 311 of \cite{Bi98ter}). Since $[\CP_{\CR}\,,\,\CP_{\CR}]\,\subset\,\h\cap\CR=\CZ$, the space $\CP_{\CR}$ is therefore an Abelian Lie sub-algebra, and, in particular Lagrangian (otherwise, its radical would be a proper $\s$-submodule in contradiction with the irreducibility of $V$).

\noindent At last, the susbspace $\CP_{\CL}$ is an $\s$-submodule of $\CP$ also isomorphic to $V$ and supplementary to $\CP_{\CR}$. Since $[\CP_{\CL},\CP_{\CL}]\subset\h\cap\CL=\s$, the subspace $\CP_{\CL}$ is Lagrangian in $\Omega$-duality with $\CP_{\CR}$. 

\noindent The action of $\CZ$ cannot be trivial. Indeed, since $[\CP_\CL\,,\,\CP_\CL]=\s$ and $[\CP_\CL\,,\,\CP_\CR]=\CZ$, Jacobi yields:
$$
[[\CP_\CL\,,\,\CP_\CL],\CP_\CR]\;\subset\;[\CZ\,,\,\CP_\CL]\;.
$$
Hence, if the right hand side is trivial, so is the action of $\s$ on $\CP_\CR\simeq V$. Simplicity of $V$ then implies item (v). 
\EPf

\noindent As an immediate corollary, one observes
\begin{cor}\label{POINC}
The symplectic manifold underlying the simply connected symplectic symmetric space $M$ associated with a generic kinematical Lie algebra of the Poincar\'e type is symplectomorphic to the cotangent bundle of a semisimple symmetric space:
$$
M\;\simeq\;T^\star(\widetilde{\L/S})
$$
where $\widetilde{\L/S}$ denotes the simply connected  symmetric space associated with the iLa $\CL=\s\oplus\CP_\CL$.

\noindent Denoting by $\L$ the connected simply connected Lie group admitting $\CL$ as a Lie algebra, the symplectomorphism can be chosen to be $\L$-equivariant.
\end{cor}
\begin{rmk}
\noindent Recall that semisimple symmetric spaces were first classified by Berger in \cite{Berger1957}.
\end{rmk}
\noindent We end this section by noticing that Theorem \ref{SSSZSS} and Corollary \ref{POINC} classify the non-flat generic kinematical Lie algebras $(\g,\s,V)$ such that $\s$ is compact semisimple:
\begin{prop}\label{CSS}
Let $(\g,\s,V)$ be a generic kinematical Lie algebras such that $\s$ is semisimple and compact. Assume that the associated simply connected symplectic symmetric space $M$ is non-flat and indecomposable. Then, $M$ is either a Hermitian symmetric space or the cotangent bundle of a Riemannian symmetric space.
\end{prop}

\begin{ex}
In order to illustrate Proposition \ref{CSS}, we note that the Poincar\'e kinematical Lie algebra is the prototypical example of such generic kinematical Lie algebras. In this case, $\s$ is the rotation Lie algebra $\s\;=\;\mathfrak{so}(D)$ naturally represented on $V=\R^D\simeq\CP_CL\simeq\CP_\CR$. The Levi factor $\CL$ is the Lorentz Lie algebra with Cartan decomposition:
$$
\CL\;=\;\s\;\oplus\;\CP_{\CL}\;.
$$
Under matrix form, $\mathfrak{so}(D)$ is represented by the lower diagonal antisymmetric $D\times D$ matrices while $\CP_\CL$ is represented by the symmetric matrices:
$$
\left(
\begin{array}{cc}
0&{}^\tau\CP_\CL\\
\CP_\CL&\mathfrak{so}(D)
\end{array}
\right)\;=\;\CL\;.
$$
\noindent The Minkowski space $\fM^{D+1}$ on which $\CL$ naturally acts is represented by the column vectors
$$
\fM^{D+1}\;=\;
\left(
\begin{array}{c}
\CZ\\
\CP_\CR
\end{array}
\right)\;.
$$
It is, to our opinion, remarkable that the Poincar\'e group therefore turns out to be the transvection group of a symplectic symmetric space. To our knowledge, this fact has not been observed in the literature yet. Perhaps partly because the Loos connection is here purely symplectic: it leaves parallel the symplectic form but it is not the Levi-Civita connection for any pseudo-Riemannian metric on the space.

\noindent It is worth noticing that the underlying symplectic manifold of our above defined symplectic symmetric space is symplectomorphic to the cotangent bundle  $T^\star(Q)$ of the hyperbolic Riemannian symmetric space $Q:=SO_o(1,D)/SO(D)$ (i.e. the mass-shell hyperboloid). 

\end{ex}

\section{Remarks on coadjoint orbits}\label{HAM}

\noindent As mentioned in the introduction,  since it is a homogeneous symplectic   
manifold, any simply connected symplectic symmetric space is a symplectic 
covering space of a coadjoint orbit of some Lie group. One may therefore wonder
which coadjoint orbits arise that way, and if they do, under what conditions 
are they symplectic symmetric spaces?
Generally, there exists no simple criterion, either geometric or algebraic, 
to decide whether a coadjoint orbit $\CO$ of a Lie group $\L$ admits a 
$\L$-invariant (local or global) symplectic symmetric space structure.

\noindent For instance, the long standing conjecture (1930) of E. Cartan stating 
that every complex homogeneous bounded domain is a symmetric space 
is a particular instance of this question in the K\"ahler case. 
In 1955, the assertion of the conjecture was proven by Borel and Matsushima 
in the case where the group of biholomorphic transformations is semisimple.
Five years later, Piateskii-Shapiro eventually answered the conjecture by 
the negative, providing a complete description of the fine structure of homogeneous 
bounded complex domains (normal j-algebras). 

\noindent More generally, proving that a coadjoint orbit $\CO$ of $\L$ admits an $\L$-invariant local symplectic symmetric space structure implies constructing an involutive Lie algebra $(\CG=\CK\oplus\CP,\sigma)$ containing $\CL$ (the Lie algebra of $\L$) and such that, denoting
$$
\pi:\CG\to\CP
$$ 
the natural projection parallel to $\CK\,$, one has
\begin{enumerate}
\item[(1)] the kernel of the restriction to $\CL$ of $\pi$ equals the Lie algebra of the stabilizer $\CC$ in $\L$ of a point $o\in\CO$.
\item[(2)] The action of $\CK$ on $\CP$ is symplectic under the identification between $\CP$ and $T_{o}(\CO)$ induced by condition (1).
\end{enumerate}
In the case $\CC$ is trivial (i.e. $\L$ is a Fr\"obenius Lie group), this task 
is currently out of reach: while the K\"ahler case follows from 
Dorfmeister-Nakajima's proof 
of the long standing Gindikin-Vinberg conjecture 
\cite{dorfmeister1988fundamental},
the general case of symplectic Lie groups \cite{baues2016symplectic} 
(i.e. Lie groups carrying a left-invariant symplectic structure) is 
open.

\vspace{2mm}

\noindent We now pass to some structure results on the coadjoint orbits involved in the setting of generic kinematical Lie algebra.
We start by recalling elementary preliminaries 
(see e.g. \cite{wallach2018symplectic}) 
on symplectic Lie group actions and on symplectic symmetric spaces. 
First, remind
\begin{dfn}
Let $G$ be a Lie group that acts on the left on a symplectic manifold $(M,\omega)$ by symplectomorphisms:
$$
\tau:G\times M\to M:(g,x)\mapsto\tau_{g}(x)\quad\quad\tau_{g}^{\star}\omega\;=\;\omega\;.
$$
Denoting by $\g$ the Lie algebra of $G$, one says that the action is {\bf weakly Hamiltonian} if every fundamental vector field of the action is Hamiltonian in the sense that for every element $X$ of the Lie algebra $\g$, there exists a smooth function $\lambda_{X}\in C^{\infty}(M)$ on $M$ such that
$$
{\rm d}\lambda_{X}\;=\;\iota_{X^{\star}}\omega
$$
where $\iota$ denotes the interior product and where $X^{\star}\in\Gamma^{\infty}(T(M))$ denotes the fundamental vector field on $M$ (``infinitesimal action''):
$$
X^{\star}_{x}\;:=\;\ddto\tau_{\exp(-t\,X)}(x)\;.
$$
When the (dual) moment map: $$\lambda:\g\to C^{\infty}(M)$$ is a homomorphism of Lie algebras (w.r.t. the Poisson bracket on $M$ associated with the symplectic structure $\omega$), one says that the action is {\bf strongly Hamiltonian} (or sometimes simply ``Hamiltonian'').
\end{dfn}
\begin{rmk}
Remind that lifting to universal covers of $G$ and $M$ and centrally extending by the Chevalley 2-cocycle of $\g$ (valued in the trivial representation of $\g$ on $\R$) associated to the dual moment map yields a systematic procedure for associating a strongly Hamiltonian action to any symplectic action.
\end{rmk}
\noindent The following result (due to Kostant) shows the relevance of co-adjoint orbits in the class of homogeneous symplectic spaces:
\begin{prop}
When strongly Hamiltonian every $G$-homogeneous symplectic space $(M,\omega)$ is a $G$-equivariant symplectic cover of a coadjoint orbit $\CO$ of $G$. The covering map is in this case the (geometrical) moment map:
$$
\CJ:M\to\CO\;\subset\;\g^{\star}:x\mapsto[X\to\lambda_{X}(x)]
$$
where $\CO$ is the coadjoint orbit of the element
$$
\xi_{o}:\g\to\R:X\mapsto\lambda_{X}(o)
$$
for any choice of base point $o$ in $M$.
\end{prop}
In that context, one has \cite{BiThese, BCG2}
\begin{prop}
Let $(M,\omega,s)$ be a simply connected symplectic symmetric space and 
$(\g=\h\oplus\CP,\sigma,\Omega)$ its corresponding transvection siLa. 
Denote by $G$ the connected simply connected Lie group admitting 
$\g$ as Lie algebra, and define the element $\widehat{\Omega}$ as the natural extension of 
$\Omega$ on $\CP$ to $\g$ by zero on $\h$:
$$
\widehat{\Omega}\;:=\;0_{\h\wedge\h}\;\oplus\;0_{\h\wedge\CP}\;\oplus\;\Omega_{\CP\wedge\CP}\;.
$$
Then:
\begin{enumerate}
\item[(i)] The element $\widehat{\Omega}$ is a Chevalley 2-cocycle of $\g$ 
associated with the trivial representation of $\g$ on $\R$.
\item[(ii)] The action of $G$ on $(M,\omega)$ is strongly hamiltonian if and only if 
$\widehat{\Omega}$ is a Chevalley coboundary, i.e., 
there exists an element $\xi_{o}$ in the dual $\g^{\star}$ of $\g$ such that
$$
\widehat{\Omega}\;=\;\delta\xi_{o}\;.
$$
In that case the associated moment map $\CJ:M\to\g^{\star}$ covers the co-adjoint orbit $\CO$ of the element $\xi_{o}$.
\item[(iii)] Assume the action of $G$ on $(M,\omega)$ to be strongly Hamiltonian. 
Denoting by $\Ad^{\flat}$ the coadjoint action, let us consider the (normal) 
group of ineffectiveness:
$$
N\;:=\;\{g\in G\;|\;\Ad^{\flat}_{g}(\xi)\;=\;\xi\;\;\forall\xi\in\CO\}\;.
$$
Denoting by $G_{\xi_{o}}$ the stabilizer of $\xi_{o}$ in $G$, the coadjoint orbit $\CO$ naturally is a symplectic symmetric space if and only if
$G_{\xi_{o}}/N$ is contained in the subgroup of $G/N$ of fixed elements under the natural involution of $G/N$ induced by $\tilde{\sigma}$ (c.f. Theorem \ref{TG}).
\end{enumerate}
\end{prop}
\noindent Back to our kinematical Lie algebras, we first note
\begin{prop}\label{KINETRANS}
Let $(\g,\s,V)$ be a generic kinematical Lie algebra such that the associated siLa $(\g,\sigma,\Omega)$ is indecomposable and non-flat. Then the transvection Lie algebra $\underline{\g}$ underlying the transvection siLa associated to $(\g,\sigma,\Omega)$ naturally carries a structure of generic kinematical Lie algebra 
$(\underline{\g},\underline{\s},V)$.
\end{prop}
\Pf First, due to the fact that $V$ is faithful, the action of $\h$ on $\CP$ is effective. Now, the Lie subalgebra $\underline{\h}:=[\CP,\CP]$ of $\h$ is necessarily transverse to $\h$ as $\Omega$ has trivial radical.
Now, consider the map $\zeta:\underline{\h}\to\CZ$ defined as the restriction to $\underline{\h}$ of the 
projection $\pr_\CZ:\h\to\CZ$ parallel to $\s$. This map is not identically zero as symplectic. Therefore
its kernel $\underline{\s}:=\ker(\zeta)$ is a co-dimension one ideal of $\underline{\h}$ contained in $\s$.
Hence $\underline{\s}=\pr_\s(\underline{\h})$ where $\pr_\s$ denotes the projection of $\h$ onto $\s$ parallel to 
$\CZ$. Which entails that $\CZ$ is contained in $\underline{\h}$.
\EPf
\begin{rmk}
Observe that if $\s$ is simple then $\g=\underline{\g}$, i.e., the siLa associated to our generic kinematical Lie algebra is, under the conditions of Proposition \ref{KINETRANS}, necessarily transvection.
\end{rmk}

\begin{prop} Let $(\g,\s,V)$ be a generic kinematical Lie algebra on the base field $\K=\R$ or $\C$. Assume the associated siLa
$(\g,\sigma,\Omega)$ to be non-flat. Denote by $\underline{G}$ be the connected simply connected Lie group admitting the transvection algebra $\underline{\g}$ as Lie algebra (c.f. Proposition\ref{KINETRANS}). Then the action of $\underline{G}$ on the simply connected symplectic symmetric space $(M,\omega,s)$ associated with the transvection siLa $(\g,\sigma,\Omega)$ is strongly Hamiltonian.
\end{prop}
\Pf The Chevalley co-boundary of the projection $\xi_o:=\pr_\CZ:\underline{\g}\to\CZ$ parallel to $\underline{\s}\oplus\CP$ consists in the extension by zero on $\underline{\h}\times\underline{\h}$ of $\Omega$.
\EPf

\begin{rmk}
\noindent As the considerations of the present work are purely algebraic, we differ the study of the group of ineffectiveness and more generally the geometrical study
of the (non-simply connected) symplectic symmetric spaces (and space-times) associated with our kinematical siLa's to a further work.
\end{rmk}

\section{Conclusions and further perspectives}\label{CONCL}

\noindent In this work, we defined the notion of generic kinematical Lie algebra and proved that
the structure of such a Lie algebra is entirely coded by a symplectic symmetric space, which we precisely 
described.

\noindent However, the overlap with the usual notion of kinematical Lie algebra as introduced in the 
physics literature is complete only in the case of space dimension $D\geq4$. In lower dimension, 
there are cases that are not covered as a consequence of the fact that item 3(a) in Definition \ref{KLA} does not hold for $\s=\mathfrak{so}(D)$ with $D\leq3$. This therefore fully  justifies to consider the analogous notion where either the $\s$-module $V$ is isomorphic to the adjoint representation of $\s$ or $\s$ is Abelian. 
We will treat those cases in a future work.

\noindent As a matter of fact, generic kinematical Lie algebras are deformations 
of one another. However, the underlying symmetric spaces governing their structure 
allow to make precise the way they are deformed into each other.
Namely, we are planning to study what are called the
\emph{homotopic structure varieties} associated by W. Bertram 
\cite{Bertram,BeBi1,BeBi2} to complex or para-complex symmetric spaces in 
the context of Jordan algebra theory.

\noindent Some of our short- and mid-term further objectives are

\begin{itemize}

\item Analyze the $\CZ$-nilpotent case without Levi condition.

\item Describe the geometry of the generalized space-times associated with our generic kinematical Lie algebras (i.e. the quotient spaces analogue to the kinematical space-times in the classical context).

\item Describe the structure of those generalized space-times associated with flat symplectic symmetric spaces.

\item Investigate the harmonic analysis attached to the natural quantization of our 
symplectic symmetric spaces, along the lines of \cite{bieliavsky2015deformation,bieliavsky2002strict}, and, in particular, 
the role of the generalized space-times in this harmonic 
analytical context.

\item Study the naturally induced geometric actions in the sense of 
Souriau \cite{Souriau1970}.

\item Investigate possible infinite dimensional generalizations.

\end{itemize}

\section{Appendix}

\subsection{Elementary representation theory}

\begin{lem}\label{DECOMP}
Let $\CG$ be a finite dimensional real Lie algebra and $\CM$ be a semi-simple real $\CG$-module. Let 
$$
\CM\;:=\;\bigoplus_{i=1}^r\;\CM_i\;=\;\bigoplus_{j=1}^{r'}\;\CM'_j
$$
be two decompostions into simple $\CG$-sub-modules. We will reserve the index $i$, respectively $j$, for the first decomposition $\CM\;=\;\bigoplus_{i=1}^r\;\CM_i$, respectively for the second one $\CM\;=\;\bigoplus_{j=1}^{r'}\;\CM'_j$. Then, 
$$
r\;=\;r'
$$
and there exists a permutation $\sigma\in S_r$ and a $\CG$-intertwiner $A\in GL(\CM)$ such that for every $i\in\{1,...,r\}$:
$$
A(\CM_i)\;=\;\CM'_{\sigma(i)}\;.
$$
\end{lem}
\Pf
For all $i\in\{1,...,r\}$ and $j\in\{1,...,r'\}$, let 
$$
\pi_{ij}:\CM_i\to\CM'_j
$$
the restriction of the projection from $\CM$ onto $\CM'_j$ 
parallel to $\oplus_{k\neq j}\CM'_k$. The maps $\pi_{ij}$ are $\CG$-intertwiners. Indeed, for every $X\in\CG$ and $m\in\CM$, one has, within the second (${}^\prime$) decomposition:
$$
[X,m]\;=\;[X,\sum_{j=1}^{r'}m'_j]\;=\;\sum_{j=1}^{r'}[X,m]_j\;.
$$
Hence one must have $[X,m'_j]=[X,m]_j$, proving the assertion.

\vspace{2mm}

\noindent Note also that, viewing all $\CM_i$ and $\CM'_j$ as a subspaces of $\CM$, one has
that, for every $i$, there must exists $j$, say $j(i)$ such that $\pi_{ij}$ is an injection. 
Indeed, let us assume that $i_0$ is such that for every $j$, $\CK_j:=\ker(\pi_{i_0j})\neq\{0\}$. Since 
$\pi_{i_0j}$ is a $\CG$-intertwiner, the subspace $\CK_j\subset\CM_{i_0}$ is a $\CG$ submodule, hence for every $j$:
$\CK_j=\CM_{i_0}$ by simplicity. This implies that $\CM_{i_0}=\{0\}$.

\noindent Therefore, one can define a function
$$
\{1,...,r\}\to\{1,...,r'\}:i\mapsto j(i)\;.
$$
Therefore $r'\leq r$. By symmetry, $r'\leq r$, hence the equality. In particular the above map defines a permutation $\sigma$ of 
$\{1,...,r\}$.

\noindent Now, by irreducibility,  $\pi_{ij}$ is either trivial or an isomorphism of simple $\CG$-modules. Therefore, the operator
$$
A\;:=\;\oplus_{i}\;\pi_{i\sigma(i)}
$$
solves the question.
\EPf

\begin{lem}
Let $W$ be a $\CG$-module, $L$ a simple $\CG$-module and $W_{(L)}$ the weak isotypic component of $L$ in $W$.
Then, the $\s$-module $W_{(L)}$ is semisimple.
\end{lem}
\Pf Let $L_1$ and $L_2$ be two different $\s$-submodules of $W_{(L)}$ each of them being isomorphic to $L$. Then, as $L$ is simple, their intersection is trivial. Denote by $L_{12}\;:=\;L_1\oplus L_2$ and let $L_3$ be a third $\s$-submodule in $W_{(L)}$ not entirely contained in $L_{12}$. The intersection $L_{12}\cap L_3$ is a proper $\s$-submodule of $L_3$, hence trivial. Inductively, one constructs a semisimple $\s$-submodule $\fL$ of the isotypic component that is maximal in the sense that there is no $\s$-submodule in $W_{(L)}$ isomorphic to $L$ and not entirely contained in $\fL$. But the isotypic component is the sum of all submodules that are isomorphic to $L$, hence each of them needs to be contained in $\fL$. Therefore $\fL$ and $W_{(L)}$ coincide.
\EPf

\noindent Let $\CG$ be a finite dimensional real Lie algebra and $\CP$ be a completely reducible $\s$-module. 

\noindent Consider a simple $\CG$-module $V$  and denote by $\CP_{(V)}$ its isotypical component in $\CP$. Remind that the isotypical component is the sum in $\CP$ of all simple $\CG$-submodules that are isomorphic to $V$.
By the previous lemma \ref{DECOMP}, all simple direct factors of any decomposition of $\CP_{(V)}$ into simple $\CG$-modules are isomorphic to $V$. Now observe
\begin{lem}\label{INTDECOMP}
\noindent For every $\CG$-submodule $\CS$ of $\CP$, one has 
$$
\CS\,\cap\,\CP_{(V)}\;=\;\CS_{(V)}\;.
$$
\end{lem}
\Pf Remind that every submodule of a completely reducible module is completely reducible and start by decomposing $\CS$ into simple $\CG$-submodules:
$$
\CS\;=:\;\oplus_{\alpha\in\Delta}\CS_{\alpha}\;.
$$
And do the same with $\CS\,\cap\,\CP_{(V)}$:
$$
\CS\,\cap\,\CP_{(V)}\;=:\;\oplus_{\beta\in\Delta^V}\CS^V_{\beta}\;.
$$
Consider the projections:
$$
\pi_{\beta\alpha}:\CS^V_{\beta}\to\CS_{\alpha}\;.
$$
Schur's lemma tells us that each of those is either an isomorphism of $\CG$-module or identically zero. Hence by the previous lemma, we may assume the inclusion:
$$
\Delta^V\;\subset\;\Delta\;.
$$
Similarly, consider a decomposition into simple $\CG$-modules:
$$
\CP_{(V)}\;=:\;\oplus_{\gamma\in\Phi}\CP_\gamma
$$
and consider the projections
$$
\pi^\prime_{\alpha\gamma}:\CS^V_{\beta}\to\CP_\gamma\;.
$$
Since at least one of them must be nontrivial, Schur tells us that every $\CG$-module must be isomorphic to $V$.\EPf

\section{Statements and Declarations: Conflict of interest and data availability}

\noindent On behalf of all authors, the corresponding author states that there is no
conflict of interest and data sharing is not applicable to this article as no datasets were generated or analysed during the current study.


\begin{thebibliography}{BCG{\etalchar{+}}06}

\bibitem[AFO18]{andrzejewski2018}
Tomasz Andrzejewski and José~M. Figueroa-O'Farrill.
\newblock Kinematical lie algebras in $2+1$ dimensions.
\newblock {\em Journal of Mathematical Physics}, 59(6):061703, 2018.
\newblock arXiv:1802.04048.

\bibitem[BB10a]{BeBi1}
Wolfgang Bertram and Pierre Bieliavsky.
\newblock Homotopes of symmetric spaces i. construction by algebras with two involutions.
\newblock {\em arXiv preprint arXiv:1011.2923}, 2010.

\bibitem[BB10b]{BeBi2}
Wolfgang Bertram and Pierre Bieliavsky.
\newblock Homotopes of symmetric spaces ii. structure variety and classification.
\newblock {\em arXiv preprint arXiv:1011.3161}, 2010.

\bibitem[BB11]{bertelson2011affine}
M{\'e}lanie Bertelson and Pierre Bieliavsky.
\newblock Affine connections and symmetry jets.
\newblock {\em arXiv preprint arXiv:1103.2300}, 2011.

\bibitem[BC16]{baues2016symplectic}
Oliver Baues and Vicente Cort{\'e}s.
\newblock {\em Symplectic Lie groups}.
\newblock Soci{\'e}t{\'e} math{\'e}matique de France, 2016.

\bibitem[BCG95]{Bi2}
P~Bieliavsky, M~Cahen, and S~Gutt.
\newblock Symmetric symplectic manifolds and deformation quantization.
\newblock In {\em Modern Group Theoretical Methods in Physics: Proceedings of the Conference in Honour of Guy Rideau}, pages 63--73. Springer, 1995.

\bibitem[BCG97]{BCG2}
P~Bieliavsky, M~Cahen, and S~Gutt.
\newblock A class of homogeneous symplectic manifolds.
\newblock {\em Contemporary mathematics}, 203:241--256, 1997.

\bibitem[BCG{\etalchar{+}}06]{BCGS}
Pierre Bieliavsky, Michel Cahen, Simone Gutt, John Rawnsley, and Lorenz Schwachh{\"o}fer.
\newblock Symplectic connections.
\newblock {\em International Journal of Geometric Methods in Modern Physics}, 3(03):375--420, 2006.

\bibitem[Ber57]{Berger1957}
Marcel Berger.
\newblock Les espaces symétriques non compacts.
\newblock {\em Annales scientifiques de l'École Normale Supérieure, Série 3}, 74(2):85--177, 1957.

\bibitem[Ber08]{Bertram}
W.~Bertram.
\newblock Homotopes and conformal deformations of symmetric spaces.
\newblock {\em J. Lie Theory}, 18:301--333, 2008.

\bibitem[BG15]{bieliavsky2015deformation}
Pierre Bieliavsky and Victor Gayral.
\newblock {\em Deformation quantization for actions of K{\"a}hlerian Lie groups}, volume 236.
\newblock American Mathematical Society, 2015.

\bibitem[Bie95]{BiThese}
Pierre Bieliavsky.
\newblock {\em Symplectic Symmetric Spaces}.
\newblock PhD thesis, Université Libre de Bruxelles, Bruxelles, Belgium, 1995.
\newblock Doctoral thesis; advisors: Simone Gutt (advisor), Michel Cahen (co-advisor).

\bibitem[Bie98a]{Bi98ter}
Pierre Bieliavsky.
\newblock Four-dimensional simply connected symplectic symmetric spaces.
\newblock {\em Geometriae Dedicata}, 69(3):291--316, 1998.

\bibitem[Bie98b]{Bi98bis}
Pierre Bieliavsky.
\newblock Semisimple symplectic symmetric spaces.
\newblock {\em Geometriae Dedicata}, 73(3):245--273, 1998.

\bibitem[Bie02]{bieliavsky2002strict}
Pierre Bieliavsky.
\newblock Strict quantization of solvable symmetric spaces.
\newblock {\em J. Symplectic Geom.}, 1(2):269--320, 2002.

\bibitem[BLL68]{Bacry:1968zf}
H.~Bacry and J.~Levy-Leblond.
\newblock {Possible kinematics}.
\newblock {\em J. Math. Phys.}, 9:1605--1614, 1968.

\bibitem[BN86]{Bacry:1986pm}
H.~Bacry and J.~Nuyts.
\newblock {Classification of Ten-dimensional Kinematical Groups With Space Isotropy}.
\newblock {\em J. Math. Phys.}, 27:2455, 1986.

\bibitem[Bor50]{borel1950sous}
Armand Borel.
\newblock Sous-groupes compacts maximaux des groupes de Lie.
\newblock {\em S{\'e}minaire Bourbaki}, 1:271--279, 1950.

\bibitem[DN88]{dorfmeister1988fundamental}
Josef Dorfmeister and Kazufumi Nakajima.
\newblock The fundamental conjecture for homogeneous k{\"a}hler manifolds.
\newblock {\em Acta Mathematica}, 161(1):23--70, 1988.

\bibitem[DR52]{de1952reductibilite}
Georges De~Rham.
\newblock Sur la r{\'e}ductibilit{\'e} d'un espace de Riemann.
\newblock {\em Commentarii Mathematici Helvetici}, 26(1):328--344, 1952.

\bibitem[Eis66]{eisenhart1966riemannian}
Luther~Pfahler Eisenhart.
\newblock {\em Riemannian Geometry}.
\newblock Princeton University Press, Princeton, NJ, 1966.
\newblock Reprint of the original 1926 edition.

\bibitem[FK94]{faraut1994analysis}
Jacques Faraut and Adam Kor{\'a}nyi.
\newblock {\em Analysis on symmetric cones}.
\newblock Oxford university press, 1994.

\bibitem[FO17]{Figueroa-OFarrill:2017sfs}
Jos{\'e} Figueroa-O'Farrill.
\newblock {Classification of kinematical Lie algebras}.
\newblock 11 2017.

\bibitem[FO18a]{ofarrill2017higher}
José~M. Figueroa-O'Farrill.
\newblock Higher-dimensional kinematical Lie algebras via deformation theory.
\newblock {\em Journal of Mathematical Physics}, 59(6):061702, 2018.
\newblock arXiv:1711.07363.

\bibitem[FO18b]{ofarrill2017defo}
José~M. Figueroa-O'Farrill.
\newblock Kinematical Lie algebras via deformation theory.
\newblock {\em Journal of Mathematical Physics}, 59(6):061701, 2018.
\newblock arXiv:1711.06111.

\bibitem[FO22]{Figueroa-OFarrill:2022nui}
Jos{\'e} Figueroa-O'Farrill.
\newblock {Non-lorentzian spacetimes}.
\newblock {\em Differ. Geom. Appl.}, 82:101894, 2022.

\bibitem[FOP19]{Figueroa-OFarrill:2018ilb}
Jos{\'e} Figueroa-O'Farrill and Stefan Prohazka.
\newblock {Spatially isotropic homogeneous spacetimes}.
\newblock {\em JHEP}, 01:229, 2019.

\bibitem[Hel78]{Helgason1978}
Sigurdur Helgason.
\newblock {\em Differential Geometry, Lie Groups, and Symmetric Spaces}, volume~80 of {\em Pure and Applied Mathematics}.
\newblock Academic Press, New York, 1978.

\bibitem[Kan85a]{kaneyuki1985}
Soji Kaneyuki.
\newblock On classification of parahermitian symmetric spaces.
\newblock {\em Tokyo Journal of Mathematics}, 8(2):473--482, 1985.

\bibitem[Kan85b]{kaneyuki1985classification}
Soji Kaneyuki.
\newblock On classification of parahermitian symmetric spaces.
\newblock {\em Tokyo Journal of Mathematics}, 8(2):473--482, 1985.

\bibitem[KN63]{KN1}
Shoshichi Kobayashi and Katsumi Nomizu.
\newblock {\em Foundations of Differential Geometry}, volume~I of {\em Interscience Tracts in Pure and Applied Mathematics}.
\newblock Interscience Publishers / John Wiley and Sons, New York--London, 1963.
\newblock MR 0152974.

\bibitem[KN96]{KN2}
Shoshichi Kobayashi and Katsumi Nomizu.
\newblock {\em Foundations of differential geometry, volume 2}, volume~2.
\newblock John Wiley \& Sons, 1996.

\bibitem[Loo69]{Lo}
Ottmar Loos.
\newblock {\em Symmetric spaces: General theory}, volume~1.
\newblock WA benjamin, 1969.

\bibitem[Mal98]{maldacena1998large}
Juan~M. Maldacena.
\newblock The large n limit of superconformal field theories and supergravity.
\newblock {\em Advances in Theoretical and Mathematical Physics}, 2(2):231--252, 1998.

\bibitem[Mor23]{Morand:2023emw}
Kevin Morand.
\newblock {Possible ambient kinematics}.
\newblock {\em J. Math. Phys.}, 64(11):112502, 2023.

\bibitem[Sou70]{Souriau1970}
Jean-Marie Souriau.
\newblock {\em Structure des syst{\`e}mes dynamiques}.
\newblock Dunod, Paris, 1970.
\newblock English translation: \emph{Structure of Dynamical Systems}, Birkh{\"a}user, 1997.

\bibitem[SV86]{SVH}
K~Sekigawa and L~Vanhecke.
\newblock Symplectic geodesic symmetries on k{\"a}hler manifolds.
\newblock {\em The Quarterly Journal of Mathematics}, 37(1):95--103, 1986.

\bibitem[Wal18]{wallach2018symplectic}
Nolan~R Wallach.
\newblock {\em Symplectic geometry and Fourier analysis}.
\newblock Courier Dover Publications, 2018.

\bibitem[Wei94]{W}
Alan Weinstein.
\newblock Traces and triangles in symmetric symplectic spaces.
\newblock {\em Symplectic geometry and quantization}, 179:261--270, 1994.

\bibitem[Wei08]{weinberg2008cosmology}
Steven Weinberg.
\newblock {\em Cosmology}.
\newblock Oxford University Press, Oxford, UK; New York, NY, 2008.

\bibitem[Wu63]{wu1963rham}
Hongxi Wu.
\newblock {\em On the de Rham decomposition theorem}.
\newblock PhD thesis, Massachusetts Institute of Technology, 1963.

\bibitem[Wu64]{deRhamWu}
H.~Wu.
\newblock On the de Rham decomposition theorem.
\newblock {\em Illinois Journal of Mathematics}, 8(2):291--311, Jun 1964.

\end{thebibliography}

\newcommand{\etalchar}[1]{$^{#1}$}

\end{document}